\documentclass[12pt]{article}
\usepackage{amssymb}

\begin{document}
\newcommand{\bea}{\begin{eqnarray}}
\newcommand{\ena}{\end{eqnarray}}
\newcommand{\beas}{\begin{eqnarray*}}
\newcommand{\enas}{\end{eqnarray*}}
\newcommand{\beq}{\begin{equation}}
\newcommand{\enq}{\end{equation}}
\def\qed{\hfill \mbox{\rule{0.5em}{0.5em}}}
\newcommand{\From}{From}
\newcommand{\ignore}[1]{}
\newcommand{\balpha}{\mbox{\boldmath {$\alpha$}}}
\newcommand{\transpose}{{\mbox{\scriptsize\sf T}}}
\newcommand{\E }{\makebox[2pt][l]{\bf E} {\rm I}~}
\newcommand{\Prob}{{\bf P}}
\newcommand{\etamu}{\mu}
\newcommand{\V}{S}
\newcommand{\id}{\mbox{id}}
\newtheorem{theorem}{Theorem}[section]
\newtheorem{corollary}[theorem]{Corollary}
\newtheorem{proposition}[theorem]{Proposition}
\newtheorem{lemma}[theorem]{Lemma}
\newtheorem{definition}[theorem]{Definition}
\newtheorem{exmpl}[theorem]{Example}
\newtheorem{remark}[theorem]{Remark}
\newtheorem{condition}[theorem]{Condition}
\newenvironment{example}{\begin{exmpl}\begin{upshape}} {\end{upshape}\end{exmpl}}

\title{{\bf\Large Berry Esseen Bounds for Combinatorial Central Limit
Theorems and Pattern Occurrences, using Zero and Size Biasing}
\thanks{AMS 2000 subject classifications. Primary 60F05, \ignore{60F05 Central limit and other weak theorems}
60C05 \ignore{60C05 Combinatorial probability}.} \,\,\thanks{Key
words and phrases: smoothing inequality, Stein's method,
permutations, graphs}}
\author{Larry Goldstein\\University of Southern California}

\maketitle{}


\begin{abstract}
Berry Esseen type bounds to the normal, based on zero- and
size-bias couplings, are derived using Stein's method. The zero
biasing bounds are illustrated with an application to
combinatorial central limit theorems where the random permutation
has either the uniform distribution or one which is constant over
permutations with the same cycle type and having no fixed points.
The size biasing bounds are applied to the occurrences of fixed
relatively ordered sub-sequences (such as rising sequences) in a
random permutation, and to the occurrences of patterns, extreme
values, and subgraphs on finite graphs.
\end{abstract}

\section{Introduction}
\label{introduction} Berry Esseen type bounds for normal
approximation are developed using Stein's method, based on zero
and size bias couplings. The results are applied to bound the
proximity to the normal in combinatorial central limit theorems
where the random permutation has either a uniform distribution, or
one which is constant over permutations with the same cycle type,
with no fixed points; to counting the number of occurrences of
fixed, relatively ordered sub-sequences, such as rising sequences,
in a random permutation; and to counting on finite graphs the
number of occurrences of patterns, local extremes, and subgraphs.

Stein's method ([\ref{Stein72}], [\ref{Stein86}]) uses
characterizing equations to obtain bounds on the error when
approximating distributions by a given target. For the normal
[\ref{Stein81}], $X \sim {\cal N}(\etamu,\sigma^2)$ if and only if
\bea
\label{characterize-Z} \E (X-\etamu)f(X) = \sigma^2 \E f'(X)
\ena
for all absolutely continuous $f$ for which $\E|f'(X)|<\infty$.
From such a characterizing equation, a difference or differential
equation can be set up to bound the difference between the
expectation of a test function $h$ evaluated on a given variable
$Y$, and on the variable $X$ having the target distribution. For
the normal, with $X$ having the same mean $\mu$ and variance
$\sigma^2$ as $Y$, the characterizing equation
(\ref{characterize-Z}) leads to the differential equation \beq
\label{basic} h((y-\etamu)/\sigma) - Nh= \sigma^2
f'(y)-(y-\etamu)f(y), \enq where $Nh=\E h(Z)$ with $Z \sim {\cal
N}(0,1)$, the standard normal mean of the test function $h$. At
$Y$, the expectation of the left hand side can be evaluated by
calculating the expectation of the right hand side using the
bounded solution $f$ of (\ref{basic}) for the given $h$. By this
device, Stein's method can handle various kinds of dependence
through the use of coupling constructions.

We consider and compare two couplings of a given $Y$ to achieve
normal bounds. First, for $Y$ with mean zero and variance
$\sigma^2 \in (0,\infty)$, we say that $Y^*$ has the $Y$-zero
biased distribution if \beq \label{zero-bias} \E Yf(Y)=\sigma^2 \E
f'(Y^*) \enq for all absolutely continuous functions $f$ for which
the expectation of either side exists. This `zero bias
transformation' from $Y$ to $Y^*$ was introduced in
[\ref{Goldstein-Reinert}], and it was shown there that $Y^*$
exists for every mean zero $Y$ with finite variance. Similarly,
for $Y$ non-negative with finite mean $\E Y= \etamu$, we say that
$Y^s$ has the $Y$-size biased distribution if
\bea
\label{size-bias} \E Yf(Y) = \etamu \E f(Y^s)
\ena
for all $f$ for which the expectation of either side exists. The
size biased distribution exists for any non-negative $Y$ with
finite mean, and was used for normal approximation in
[\ref{Goldstein-Rinott}].

A coupling $(Y,Y^*)$ where $Y^*$ has the $Y$-zero biased
distribution lends itself for use in the Stein equation
(\ref{basic}) in the following way; by (\ref{zero-bias}), with
$\sigma^2=1$ say, we have
\bea
\label{W-to-W*-via-f} \E h(Y) - Nh = \E \left[ f'(Y) - Yf(Y)
\right] = \E \left[ f'(Y) - f'(Y^*)\right].
\ena
Therefore, the difference between $Y$ and the normal, as tested on
$h$, equals the difference between $Y$ and $Y^*$, as tested on
$f'$. Additionally, as observed in [\ref{Goldstein-Reinert}] and
seen directly from (\ref{W-to-W*-via-f}), $Y$ is normal if and
only if $Y=_d Y^*$. It is therefore natural that the distance from
$Y$ to the normal can be expressed in terms of distance from $Y$
to $Y^*$. Theorem \ref{delta-bound-theorem} makes this statement
precise, showing that the distance from the standardized $Y$ to
the normal as measured by $\delta$ in (\ref{def-delta}) depends on
the distribution of $Y$ only through a bound on $|Y^*-Y |$. A
similar phenomenon is seen in [\ref{res}] with $d_{{\rm W}}$ the
Wasserstein distance, where it is shown that, for any mean zero
variance $\sigma^2$ variable $Y$, and $X \sim {\cal
N}(0,\sigma^2)$,

\beas
d_{{\rm W}}(Y,X) \le 2 d_{{\rm W}}(Y,Y^*).
\enas

The use of size bias couplings in the Stein equation in
(\ref{size-bias-stein-equation-connection}),
(\ref{four-terms-in-two-groups}) and subsequent calculations
depends on the following identity, which is applied in a less
direct manner than (\ref{W-to-W*-via-f}); for $Y \ge 0$ with mean
$\mu$ and variance $\sigma^2$,
$$
\E(Y-\mu)f(Y)=\mu \E \left(f(Y^s)-f(Y)\right) \quad \mbox{and
therefore} \quad \sigma^2 = \mu \E(Y^s-Y).
$$

With $W=(Y-\etamu)/\sigma$, many authors (e.g. [\ref{Bolt}],
[\ref{Gotze}], [\ref{BG93}], [\ref{Rin-Rot-1}], [\ref{Rin-Rot-2}],
[\ref{Chen-Shao}]) have been successful in obtaining bounds on the
distance
\bea
\label{def-delta} \delta = \sup_{h \in {\cal H}}|\E h(W) - Nh|
\ena
to the normal, over classes of non-smooth functions ${\cal H}$,
using Stein's method. Here we take the smoothing inequality
approach, following [\ref{Rin-Rot-2}]. In particular, ${\cal H}$
is a class of measurable functions on the real line such that

\begin{enumerate}
\item[(i)] The  functions $h \in {\cal H }$ are uniformly bounded in
absolute value by a constant, which we take to be 1 without loss
of generality,

\item[(ii)]  For any real numbers $c$ and $d$, and for any $h(x)
\in {\cal H }$, the function $h(cx+d) \in {\cal H }$,

\item[(iii)] For any $\epsilon >0$ and $h \in {\cal H }$, the
functions $h^+_\epsilon, \,  h^-_\epsilon$ are also in ${\cal H
}$, and
\begin{equation}
\label{aaa} \E \tilde{h}_\epsilon(Z) \leq a\epsilon
\end{equation}
for some constant $a$ which depends only on the class ${\cal H
}$, where
\bea
\label{H-closure} h^+_{\epsilon}(x) = \sup_{|y| \leq
\epsilon}h(x+y), \quad h^-_{\epsilon}(x) = \inf_{|y| \leq
\epsilon}h(x+y), \quad \mbox{and} \quad \tilde{h}_{\epsilon}(x) =
h^+_{\epsilon}(x)-h^-_{\epsilon}(x).
\ena
\end{enumerate}

The collection of indicators of all half lines, and indicators of
all intervals, for example, each form classes ${\cal H}$ which
satisfy (\ref{H-closure}) and (\ref{aaa}) with $a =
\sqrt{2/\pi}$ and $a=2\sqrt{2/\pi}$ respectively (see e.g.
[\ref{Rin-Rot-2}]).

Since the bound on $\delta$ in Theorem \ref{delta-bound-theorem}
depends only the size of $|Y^*-Y|$, it may be computed without the
need for the calculation of the variances of certain conditional
expectations that arise in other versions of Stein's method,
$\sqrt{\mbox{Var}\{ \E[(Y''-Y')^2|Y']\} }$ for the exchangeable
pair method, or the term (\ref{def-Delta}) for the size bias
coupling studied here.

\begin{theorem}
\label{delta-bound-theorem} Let $Y$ be a mean zero random variable
with variance $\sigma^2 \in (0,\infty)$, and $Y^*$ be defined on
the same space having the $Y$-zero biased distribution. If
$|Y^*-Y| \le 2B$ for some $B \le \sigma/24$, then for $\delta$ as
in (\ref{def-delta}) and $a$ as in (\ref{aaa}),
\bea
\label{delta-bound} \delta \le  A \left(37 + 12 A + 112a\right),
\ena
for $A=2B/\sigma$. For indicators of all half lines, and the
indicators of all intervals, using $a = \sqrt{2/\pi}$ and
$a=2\sqrt{2/\pi}$, we have respectively
\bea
\label{delta-bound-cases}
\delta \le  A \left(127 + 12 A \right)
\quad \mbox{and} \quad \delta \le  A \left(216 + 12 A \right).
\ena
\end{theorem}
See (\ref{delta-bound-alt}) and (\ref{size-delta-bound-alt}) for
some variations on the bound (\ref{delta-bound}) here, and
(\ref{size-delta-bound}) below, respectively. We note that Theorem
\ref{delta-bound-theorem} immediately provides a bound on $\delta$
of order $\sigma^{-1}$ whenever $|Y^*-Y|$ is bounded. In Section
\ref{zero-applications}, we apply Theorem
\ref{delta-bound-theorem} to random variables of the form
\bea
\label{comb-W}
Y=\sum_{i=1}^n a_{i,\pi(i)},
\ena
depending on a fixed array of real numbers $\{a_{ij}\}_{i,j=1}^n$
and a random permutation $\pi \in {\cal S}_n$, the symmetric
group. In Section \ref{comb-clt-uniform} we consider $\pi$ having
the uniform distribution on ${\cal S}_n$, and in Section
\ref{sec-cycle} distributions constant on cycle type having no
fixed points (conditions (\ref{constant-on-conjugacy}) and
(\ref{no-one-cycles}) respectively).

For a size bias coupling $(Y,Y^s)$, Theorem
\ref{size-delta-bound-theorem} gives a bound on $\delta$ which
depends on the size of $|Y^s-Y|$, and additionally on $\Delta$ in
(\ref{def-Delta}). While $\Delta$ may be difficult to calculate
precisely in many cases, size bias couplings can be more easily
constructed for a broader range of examples than the zero bias
couplings.

\begin{theorem}
\label{size-delta-bound-theorem} Let $Y \ge 0$ be a random
variable with mean $\etamu$ and variance $\sigma^2  \in
(0,\infty)$, and let $Y^s$ be defined on the same space, with the
$Y$-size biased distribution. If $|Y^s-Y| \le B$ for some $B \le
\sigma^{3/2}/\sqrt{6 \etamu}$, then for $\delta$ as in
(\ref{def-delta}) and $a$ as in (\ref{aaa}),
\bea
\label{size-delta-bound} \delta \le \frac{aA}{2} +
\frac{\etamu}{\sigma} \left((19+56a)A^2 + 4A^3 \right)+\frac{23
\mu \Delta }{\sigma^2},
\ena
where
\bea
\label{def-Delta}
\Delta = \sqrt{\mbox{Var}(\E (Y^s-Y|Y))}
\ena
and $A=B/\sigma$. For indicator functions of all half lines and
the indicators functions of all intervals, by using $a =
\sqrt{2/\pi}$ and $a=2\sqrt{2/\pi}$, we respectively find that
\beas
\delta \le 0.4A +
\frac{\mu}{\sigma}\left(64 A^2+4A^3\right) + \frac{23 \Delta
\mu}{\sigma^2} \quad \mbox{and} \quad \delta \le  0.8 A +
\frac{\mu}{\sigma}\left(109 A^2+4A^3\right) + \frac{23 \Delta
\mu}{\sigma^2}.
\enas
\end{theorem}
If the mean $\etamu$ is of order $\sigma^2$, $B$ is bounded and
$\Delta=\sigma^{-1}$, then $\delta$ will have order
$O(\sigma^{-1})$. The application of Theorem
\ref{size-delta-bound-theorem} to counting the occurrences of
fixed relatively ordered sub-sequences, such as rising sequences,
in a random permutation, and to counting the occurrences of color
patterns, local maxima, and sub-graphs  in finite graphs is
illustrated in Section \ref{size-applications}. The proofs of
Theorems \ref{size-delta-bound-theorem} and
\ref{delta-bound-theorem} are given in Section \ref{proofs}.

Nothing should be inferred from the fact that the zero bias
applications presented here involve global dependence, and that
the dependence in the examples used to illustrate the size bias
approach is local; the exchangeable pair coupling on which our
zero biased constructions are based can also be applied in cases
of local dependence, and the size bias approach was applied in
[\ref{Goldstein-Rinott}] to variables having global dependence.

In both zero and size biasing, a sum $Y=\sum_{\alpha \in {\cal A}}
X_\alpha$ of independent variables on a finite index set ${\cal
A}$ is biased by choosing a summand at random and replacing it
with its biased version. To describe the zero biasing coupling,
let $\{X_\alpha\}_{\alpha \in {\cal A}}$ be a collection of mean
zero variables with finite variance, and $I$ an independent random
index with distribution
\bea
\label{Ipropvar}
P(I=\alpha) = \frac{w_\alpha}{\sum_{\beta \in {\cal A}}w_\beta},
\ena
where $w_\alpha=\mbox{Var}(X_\alpha)$. It was shown in
[\ref{Goldstein-Reinert}] that replacing $X_I$ by a variable
$X_I^*$ having the $X_I$-zero bias distribution, independent of
$\{X_\alpha, \alpha \not = I\}$, gives
\bea
\label{W*=W-X+X}
Y^*=Y-X_I+X_I^*,
\ena
a variable having the $Y$-zero biased distribution. Hence, when a
sum of many independent variables of the same order is coupled
this way to its zero biased version, the magnitude of
$(Y^*-Y)/\sqrt{\mbox{Var}(Y)}$, and therefore of distance measures
such as $\delta$, are small.

The construction of the size biased coupling in the independent
case is similar. Let $\{X_\alpha\}_{\alpha \in {\cal A}}$ be a
collection of non-negative variables with finite mean. Then, with
$I$ a random index independent of all others variables, having
distribution (\ref{Ipropvar}) with $w_\alpha=\E X_\alpha$, the
replacement of $X_I$ by a variable $X_I^s$ with the $X_I$-size
bias distribution, independent of the remaining variables, gives a
variable with the $Y$-size biased distribution.

Zero biased couplings of $Y^*$ to a sum $Y$ of non-independent
variables $X_1,\ldots,X_n$ is presently not very well understood.
A construction in the presence of the weak global dependence of
simple random sampling was given in [\ref{Goldstein-Reinert}].
Based on a remark in [\ref{Goldstein-Reinert}], we here exploit a
connection between the zero bias coupling and the exchangeable
pair $(Y',Y'')$ of [\ref{Stein86}] with distribution $dP(y',y'')$
satisfying $\E(Y''|Y')=(1-\lambda)Y'$ for some $\lambda \in
(0,1]$; in particular, we make use of a pair
$(Y^\dagger,Y^\ddagger)$ with distribution proportional to
$(y'-y'')^2 dP(y',y'')$.

The construction of $Y$ and $Y^s$ on a common space for the sum of
non-independent variables $X_1,\ldots,X_n$ is more direct, and was
described in Lemma 2.1 of [\ref{Goldstein-Rinott}]; we choose a
summand with probability proportional to its expectation, replace
it by one from its size-biased distribution, and then adjust the
remaining variables according to the conditional distribution
given the value of the newly chosen variable. This construction is
applied in Section \ref{size-applications}, and a `squared' zero
biasing form of it in Section \ref{zero-applications}.

The mappings of a distribution $Y$ to its zero biased $Y^*$ or
size biased $Y^s$ versions are special cases of distributional
transformations from $Y$ to some $Y^{(m)}$ which are specified by
a function $H$ and characterizing equation
$$
\E H(Y)f(Y) = \eta \E f^{(m)}(Y^{(m)}) \quad \mbox{for all smooth
$f$,}
$$
where $f^{(m)}$ denotes the $m^{th}$ derivative of $f$, and $\eta$
is, necessarily, $(m!)^{-1}\E H(Y)Y^m$ when this expectation
exists. The zero bias and size bias transformation correspond to
$m=1$ and $H(x)=x$, and $m=0$ and $H(x)=x^+$, respectively. In
general, such a $Y^{(m)}$ exists when $H$ and $Y$ satisfy certain
sign change and orthogonality properties, as discussed in
[\ref{Goldstein-Reinert-03}].

\section{Zero Biasing: Combinatorial Central Limit Theorems}
\label{zero-applications}

In this section, we illustrate the use of Theorem
\ref{delta-bound-theorem} to obtain Berry Esseen bounds in
combinatorial central limit theorems, that is, for variables $Y$
as in (\ref{comb-W}), in Section \ref{comb-clt-uniform} we do so
for permutations having the uniform distribution over the
symmetric group and, in Section \ref{sec-cycle}, we do so for
permutations with distribution constant on those having the same
cycle type, with no fixed points. First we present Proposition
\ref{peara}, which suggests a method for the construction of zero
bias couplings based on the existence of exchangeable pairs; its
statement appears in [\ref{Goldstein-Reinert}].
\begin{proposition}
\label{peara} Let $Y'$ and $Y''$ be an exchangeable pair, with
distribution $dP(y',y'')$ and $\mbox{Var}(Y')=\sigma^2 \in
(0,\infty)$, which satisfies the linearity condition
\bea
\label{linear-lambda-F} \E(Y''|Y')=(1-\lambda)Y' \quad \mbox{for
some $\lambda \in (0,1]$.}
\ena
Then
\bea
\label{simple-moment-claims} \E Y'=0 \quad \mbox{and} \quad \E
(Y'-Y'')^2=2\lambda \sigma^2,
\ena
and if $Y^\dagger$ and $Y^\ddagger$ have distribution
\bea
\label{Z-distribution} dP^\dagger(y',y'') = \frac{(y'-y'')^2 }{\E
(Y'-Y'')^2}dP(y',y''),
\ena
and $U \sim {\cal U}[0,1]$ is independent of $Y^\dagger$ and
$Y^\ddagger$, then the variable
\bea
\label{UZ+(1-U)Z} Y^* = U Y^\dagger + (1-U)Y^\ddagger \quad
\mbox{has the $Y'$ zero biased distribution.}
\ena
\end{proposition}
\noindent {\bf Proof:} The claims in (\ref{simple-moment-claims})
follow from (\ref{linear-lambda-F}) and exchangeability. Hence we
need only show that $Y^*$ in (\ref{UZ+(1-U)Z}) satisfies
(\ref{zero-bias}). For a differentiable test function $f$,
\beas
\sigma^2 \E f'(U Y^\dagger + (1-U)Y^\ddagger) &=& \sigma^2 \E
\left( \frac{f(Y^\dagger)-
f(Y^\ddagger)}{Y^\dagger-Y^\ddagger}\right) \\ &=& \sigma^2 \E
\left( \frac{(Y'-Y'')(f(Y')- f(Y''))}{\E (Y'-Y'')^2}\right).
\enas
Now if we use (\ref{linear-lambda-F}) to obtain $\E
Y''f(Y')=(1-\lambda)\E Y'f(Y')$, followed by
(\ref{simple-moment-claims}), expanding yields
\beas
\sigma^2 \E  \left( \frac{Y'f(Y')-Y''f(Y')-Y'f(Y'')+Y''f(Y'')}{\E
(Y'-Y'')^2} \right)=\frac{2\lambda \sigma^2 \E Y'f(Y')}{\E
(Y'-Y'')^2}=\E Y'f(Y'). \qed
\enas

\begin{example}
\label{Y'Y''independent} Given a mean zero finite variance $Y'$,
let $Y''$ be an independent copy of $Y'$. The pair $(Y',Y'')$
satisfies the conditions of Proposition \ref{peara} with
$\lambda=1$, and hence, $Y^*$ as in (\ref{UZ+(1-U)Z}) has the $Y'$
zero bias distribution with $(Y^\dagger,Y^\dagger)$ as in
(\ref{Z-distribution}). However, by coupling $Y'$ close to $Y''$,
so that $Y^\dagger$ is close to $Y^\ddagger$, causes $B$, and
therefore, the bound $\delta$ of Theorem \ref{delta-bound-theorem}
to be small.
\end{example}

\begin{remark}
\label{suggested-construction}
 The following construction of
$(Y^\dagger,Y^\ddagger)$ suggested by Proposition \ref{peara} is
similar to the one used for size biasing (see Lemma 2.1 of
[\ref{Goldstein-Rinott}] and Section \ref{size-applications}).
Given $Y'$, first construct an exchangeable $Y''$ close to $Y'$
satisfying (\ref{linear-lambda-F}), and then, independently
construct the variables appearing in the `square biased' term
$(Y'-Y'')^2$. Lastly, adjust the remaining variables that make up
$(Y',Y'')$ to have their original conditional distribution, given
the newly generated variables.
\end{remark}

\begin{example}
Let $\{X_i',X_i''\}_{i=1,\ldots,n}$ be i.i.d. mean zero variables
with finite variancess, let $Y'=\sum_{i=1}^n X_i'$, and let $I$ be
an independent random index with uniform distribution over
$\{1,\ldots,n\}$. Letting $Y''=Y'-X_I'+X_I''$, the pair $(Y',Y'')$
is exchangeable and satisfies the conditions of Proposition
\ref{peara}, with $\lambda=1/n$. Set $S=\sum_{i \not = I}X_i'$ and
$(T',T'')=(X_I',X_I'')$. Applying Example \ref{Y'Y''independent}
to $(T',T'')$, and forming $(T^\dagger,T^\ddagger)$ independently
of $\{X_i',X_i''\}_{i \not = I}$, gives $UT^\dagger +
(1-U)T^\ddagger=X_I^*$. By their independence from $X_I',X_I''$,
$\{X_i',X_i''\}_{i \not = I}$ already
 have their original conditional distribution, given
$(T^\dagger, T^\ddagger)$; hence $Y^*=\V+X_I^*$, in agreement with
(\ref{W*=W-X+X}).
\end{example}

Applying this construction in the presence of dependence results
in $\V$, a function of the variables which can be kept fixed, and
variables $T',T^\dagger,T^\ddagger$, on a joint space, such that
\bea \label{new-def-T} Y'=\V+T', \quad Y^\dagger=\V+T^\dagger,
\quad \mbox{and} \quad Y^\ddagger=\V+T^\ddagger. \ena When
$T',T^\dagger$ and $T^\ddagger$ are all bounded by $B$,
(\ref{UZ+(1-U)Z}) gives \bea \label{Y*Y2B-bound}
|Y^*-Y'|&=&|UT^\dagger+(1-U)T^\ddagger-T'| \le
U|T^\dagger|+(1-U)|T^\ddagger| + |T'| \le 2B. \ena

Let an array $\{a_{ij}\}_{i,j=1}^n$ of real numbers satisfy \bea
\label{C-is-sup} \sum_{j=1}^n a_{ij}=0 \quad \mbox{for all $i$,
and set} \quad C=\max_{i,j}|a_{ij}|. \ena By replacing $Y$ in
(\ref{comb-W}) by $Y-\E Y$ we assume, without loss of generality,
that $\E a_{i,\pi(i)}=0$ for every $i$. In Theorem
\ref{comb-uniform-theorem}, below, where $\pi$ is uniformly
distributed over ${\cal S}_n$, this assumption is equivalent to
(\ref{C-is-sup}). In Theorem \ref{comb-constant-theorem}, since
$\pi$ has no fixed points, by (\ref{no-one-cycles}), without loss
of generality we have $a_{ii}=0$ for all $i$ in (\ref{aij-props}).
In addition, since the distribution of $\pi$ is constant on
permutations having the same cycle type, by
(\ref{constant-on-conjugacy}), $\E a_{i,\pi(i)}=(1/(n-1))\sum_{j
\not = i} a_{ij}$, and the mean zero assumption is again
equivalent to (\ref{C-is-sup}). Avoiding trivial cases, we also
assume that $\mbox{Var}(Y)=\sigma^2>0$. For ease of notation we
write $Y'$ and $\pi'$ interchangeably for $Y$ and $\pi$,
respectively, in the remainder of this section .

In Sections \ref{comb-clt-uniform} and \ref{sec-cycle} the
construction above produces variables $Y',Y^\dagger$ and
$Y^\ddagger$, given by (\ref{comb-W}) (with $\pi$ replaced by
$\pi',\pi^\dagger$ and $\pi^\ddagger$, respectively), and a set of
indices ${\cal I}$ outside of which these permutations agree, such
that (\ref{new-def-T}) holds with \bea \label{def-T'} \V=\sum_{i
\not \in {\cal I}}a_{i,\pi'(i)}, \,\, T'=\sum_{i \in {\cal
I}}a_{i,\pi'(i)}, \,\, T^\dagger=\sum_{i \in {\cal
I}}a_{i,\pi^\dagger(i)}, \,\, \mbox{and} \,\, T^\ddagger=\sum_{i
\in {\cal I}}a_{i,\pi^\ddagger(i)}. \ena Therefore $B$ in
(\ref{Y*Y2B-bound}) can be set equal to $C$ in (\ref{C-is-sup})
times a worst case bound on the size of ${\cal I}$. The
specifications of $\pi',\pi'',\pi^\dagger,$ and $\pi^\ddagger$ are
given in terms of transpositions $\tau_{ij}$, those permutations
satisfying $\tau_{ij}(i)=j,\tau_{ij}(j)=i$ and $\tau_{ij}(k)=k$
for all $k \not \in \{i,j\}$.

\subsection{Uniform permutation distribution}
\label{comb-clt-uniform} Many authors (e.g. [\ref{Bahr}]
[\ref{Bolt}], [\ref{ChenHo}]) have considered normal approximation
to the distribution of (\ref{comb-W}) when $\pi$ is a permutation
chosen uniformly from ${\cal S}_n$. In Theorem
\ref{comb-uniform-theorem}, the dependence of $\delta$ on $C$ is
not as refined as the bound in [\ref{Bolt}], which depends on an
(unspecified) universal constant times the normalized absolute
third moments of the $\{a_{ij}\}_{i,j = 1}^n$. Here, on the other
hand, an explicit constant is provided.
\begin{theorem}
\label{comb-uniform-theorem} With $n \ge 3$, let $\{a_{ij}\}_{i,j
= 1}^n$ satisfy (\ref{C-is-sup}) and let $\pi$ be a random
permutation with uniform distribution over ${\cal S}_n$. Then,
with $C$ as in (\ref{C-is-sup}), conclusions (\ref{delta-bound})
and (\ref{delta-bound-cases}) of Theorem \ref{delta-bound-theorem}
hold for the sum $Y=\sum_{i=1}^n a_{i,\pi(i)}$ with $A =
8C/\sigma$ when $A \le 1/12$.
\end{theorem}

\noindent {\bf Proof:} Given $\pi'$, take $(I,J)$ to be
independent of $\pi'$, uniformly over all pairs with $1 \le I \not
= J \le n$, and set $\pi''=\pi'\tau_{I,J}$. In particular,
$\pi''(i)=\pi'(i)$ for $i \not \in \{I,J\}$; the variables $Y'$
and $Y''$, given by (\ref{comb-W}) with $\pi'$ and $\pi''$
respectively, are exchangeable; and \bea
\label{Y-difference-uniform} Y'-Y'' =
(a_{I,\pi'(I)}+a_{J,\pi'(J)}) - (a_{I,\pi'(J)}+a_{J,\pi'(I)}).
\ena The linearity condition (\ref{linear-lambda-F}) is satisfied
with $\lambda=2/(n-1)$, since, from (\ref{Y-difference-uniform})
and (\ref{C-is-sup}),
\beas E\left( Y'-Y''|\pi' \right) &=& 2 \left( \frac{1}{n}
\sum_{i=1}^n a_{i,\pi'(i)}- \frac{1}{n(n-1)}\sum_{i \not =
j}a_{i,\pi'(j)}
\right)\\
&=& 2 \left( \frac{1}{n} \sum_{i=1}^n a_{i,\pi'(i)} +
\frac{1}{n(n-1)}\sum_{i = 1}^n a_{i,\pi'(i)} \right) =
\frac{2}{n-1}Y'.
\enas

To construct $(Y^\dagger,Y^\ddagger)$ with distribution
proportional to $(y'-y'')^2 dP(y', y'')$, choose
$I^\dagger,K^\dagger,J^\dagger,L^\dagger$ with distribution
proportional to the squared difference $(Y'-Y'')^2$, that is,
$$
P(I^\dagger=i,K^\dagger=k,J^\dagger=j,L^\dagger=l) \sim
\left[(a_{ik}+ a_{jl}) -(a_{il}+a_{jk})\right]^2,
$$
and let
\beas
\pi^\dagger = \left\{
\begin{array}{ll}
\pi  \tau_{\pi^{-1}(K^\dagger),J^\dagger}  & \mbox{if $L^\dagger=\pi(I^\dagger),K^\dagger \not =\pi(J^\dagger)$}\\
\pi  \tau_{\pi^{-1}(L^\dagger),I^\dagger}  & \mbox{if $L^\dagger \not =\pi(I^\dagger),K=\pi(J^\dagger)$}\\
\pi  \tau_{\pi^{-1}(K^\dagger),I^\dagger} \tau_{\pi^{-1}(L^\dagger),J^\dagger} & \mbox{otherwise,}\\
\end{array}
\right. \enas and $\pi^\ddagger=\pi^\dagger
\tau_{I^\dagger,J^\dagger}$. Then (\ref{new-def-T}) and
(\ref{def-T'}) hold with ${\cal
I}=\{I^\dagger,\pi^{-1}(K^\dagger),J^\dagger,
\pi^{-1}(L^\dagger)\}$, a set of size at most 4, so by
(\ref{Y*Y2B-bound}), $|Y^*-Y| \le 8C$.$\qed$

\subsection{Permutations with distribution constant over cycle type}
\label{sec-cycle}

In this section we focus on the normal approximation of $Y$ as in
(\ref{comb-W}) when the distribution of the random permutation
$\pi$ is a function only of its cycle type. Our framework includes
the case considered in [\ref{Kolchin}], the uniform distribution
over permutations with a single cycle.

Consider a permutation $\pi \in {\cal S}_n$ represented in cycle
form; in ${\cal S}_7$ for example, $\pi=((1,3,7,5),$ $(2,6,4))$ is
the permutation consisting of one 4 cycle in which $1 \rightarrow
3 \rightarrow 7 \rightarrow 5 \rightarrow 1$ and one 3 cycle where
$2 \rightarrow 6 \rightarrow 4 \rightarrow 2$. For $q=1,\ldots,n$,
let $c_q(\pi)$ be the number of $q$ cycles of $\pi$. We say
permutations $\pi$ and $\sigma$ are of the same cycle type if
$c_q(\pi)=c_q(\sigma)$ for all $q=1,\ldots,n$; $\pi$ and $\sigma$
are of the same cycle type if and only if $\pi$ and $\sigma$ are
conjugate, i.e. if and only if there exists a permutation $\rho$
such that $\pi=\rho^{-1} \sigma \rho$. Hence, we say a probability
measure $\Prob$ on ${\cal S}_n$ is constant over cycle type if
\bea \label{constant-on-conjugacy} \Prob(\pi)=\Prob(\rho^{-1}\pi
\rho) \quad \mbox{for all $\pi,\rho \in {\cal S}_n$.} \ena

In [\ref{Goldstein-Rinott-Perm}], the authors consider a
statistical test for determining when a given pairing of $n=2m$
observations shows an unusually high level of similarity; the test
statistic is of the form (\ref{comb-W}), and, under the null
hypothesis of no distinguished pairing, the distribution $\Prob$
satisfies (\ref{constant-on-conjugacy}) with $\Prob(\pi)$ equal to
a constant if $\pi$ has $m$ 2-cycles, and $\Prob(\pi)=0$
otherwise; that is, under the null, $\Prob$ is uniform over
permutations having $m$ 2-cycles. Bounds between the normal and
the null distribution of $Y$ were determined in
[\ref{Goldstein-Rinott-Perm}] using a construction in which an
exchangeable $\pi''$ is obtained from $\pi$ by a transformation
which preserves the $m$ 2-cycle structure. The construction in
Theorem \ref{comb-constant-theorem} preserves the cycle structure
in general and, when there are $m$ 2-cycles, specializes to one
similar, but not equivalent, to that of
[\ref{Goldstein-Rinott-Perm}].

\begin{theorem}
\label{comb-constant-theorem} With $n \ge 4$, let an array
$\{a_{ij}\}_{ij = 1}^n$ of real numbers satisfy (\ref{C-is-sup}),
let
\bea
\label{aij-props} a_{ij}=a_{ji} \quad \mbox{and} \quad a_{ii}=0,
\ena
and let $\pi \in {\cal S}_n$ be a random permutation with
distribution $\Prob$ constant on cycle type, with no fixed points.
That is, $\Prob$ satisfies (\ref{constant-on-conjugacy}),
(\ref{aij-props}), and
\bea
\label{no-one-cycles} \Prob(\pi)=0 \quad \mbox{if $c_1(\pi) \not =
0$.}
\ena
Then, with $C$ as in (\ref{C-is-sup}), conclusions
(\ref{delta-bound}) and (\ref{delta-bound-cases}) of  Theorem
\ref{delta-bound-theorem} hold for the sum $Y=\sum_{i=1}^n
a_{i,\pi(i)}$ with $A = 40C/\sigma$ when $A \le 1/12$.
\end{theorem}

\noindent {\bf Proof}: To fully highlight the reason for the
imposition of the conditions (\ref{aij-props}) and
(\ref{no-one-cycles}), and also to make the complete case analysis
easier to follow, we initially consider an array satisfying only
the consequence $\sum_{1 \le i,j \le n}a_{ij}=0$ of
(\ref{C-is-sup}), and a $\Prob$ not necessarily satisfying
(\ref{no-one-cycles}).

Again, using the construction outlined in Remark
\ref{suggested-construction}, we first construct $\pi''$ from the
given $\pi'$. Let $I$ and $J$, $1 \le I \not = J \le n$ be chosen
uniformly and independently of $\pi'$, and let
$\pi''=\tau_{IJ}\pi'\tau_{IJ}$; that is, $\pi''$ is obtained by
interchanging $I$ and $J$ in the cycle representation of $\pi'$.
We claim the pair $\pi',\pi''$ is exchangeable. For fixed
permutations $\sigma'', \sigma'$, if $\sigma' \not =
\tau_{IJ}\sigma'' \tau_{IJ}$ then
$$
\Prob(\pi''=\sigma'',\pi'=\sigma')=0=\Prob(\pi'=\sigma'',\pi''=\sigma').
$$
Otherwise, $\sigma' = \tau_{IJ}\sigma'' \tau_{IJ}$ and, using
(\ref{constant-on-conjugacy}) for the second equality, we have
\beas \Prob(\pi''=\sigma'',\pi'=\sigma') = \Prob(\pi'=\sigma') =
\Prob(\pi'=\tau_{IJ}\sigma'\tau_{IJ})=\Prob(\pi''=\sigma')=\Prob(\pi'=\sigma'',\pi''=\sigma').
\enas Consequently, $Y$ and $Y''$, given by (\ref{comb-W}) with permutations
$\pi$ and $\pi''$, respectively, are exchangeable. By conditioning
on $\pi$, we show $Y',Y''$ satisfies the linearity condition
(\ref{linear-lambda-F}) with $\lambda=4/n$.

Let $S$ be the size of the set $\{I,J,\pi(I),\pi(J)\}$, and, for
$i \in \{1,\ldots,n\}$ let $|i|$ denote the number of elements in
the cycle of $\pi$ that contains $i$. Since $I \not =J$, we have
$2 \le S \le 4$. When $S=2$, either $\pi(I)=I$ and $\pi(J)=J$, or
$\pi(I)=J$ and $\pi(J)=I$; in the both cases $\pi''=\pi$. There
are four cases for $S=3$; either $A_{I,J}=\{|I|=1,|J| \ge 2\}$ or
$I$ and $J$ are interchanged (denoted by $A_{J,I}$); or $I,J$ and
$\pi(J)$ are three consecutive distinct values of $\pi$, indicated
by $B_{I,J}$, or $I$ and $J$ are interchanged (denoted by
$B_{J,I}$). The case $S=4$ is indicated by $F$. Hence,
\bea \label{A} Y'-Y''&=& \left(
a_{I,I}+a_{\pi^{-1}(J),J}+a_{J,\pi(J)} -
(a_{J,J}+a_{\pi^{-1}(J),I}+a_{I,\pi(J)}) \right)
A_{I,J}\\
\nonumber
&+&\left( a_{J,J}+a_{\pi^{-1}(I),I}+a_{I,\pi(I)} -
(a_{I,I}+a_{\pi^{-1}(I),J}+a_{J,\pi(I)}) \right)
A_{J,I}\\
\nonumber
&+&\left( a_{\pi^{-1}(I),I}+a_{I,J}+a_{J,
\pi(J)}-(a_{\pi^{-1}(I),J}+a_{J,I}+a_{I, \pi(J)})
\right)B_{I,J}\\
\nonumber
&+&\left( a_{\pi^{-1}(J),J}+a_{J,I}+a_{I,
\pi(I)}-(a_{\pi^{-1}(J),I}+a_{I,J}+a_{J, \pi(I)})
\right)B_{J,I}\\
\nonumber
&+&\left( a_{\pi^{-1}(I),I}+a_{I,\pi(I)} +
a_{\pi^{-1}(J),J}+a_{J,\pi(J)}- (a_{\pi^{-1}(I),J}+a_{J,\pi(I)}
+ a_{\pi^{-1}(J),I}+a_{I,\pi(J)}) \right) F.
\ena

For example, using the fact that the sum of $a_{\pi^{-1}(J),J}$ is
the same as that of $a_{J,\pi(J)}$ over a given cycle, the
contribution to $n(n-1)\E(Y'-Y''|\pi)$ from $A_{I,J}=\{|I|=1,|J|
\ge 2\}$, added to the equal one from $A_{J,I}$, simplifies to
\bea \label{A3terms} 2(n-3c_1(\pi))
\sum_{|i|=1}a_{i,i} + 4c_1(\pi) \sum_{|i| \ge 1} a_{i,\pi(i)} -
2c_1(\pi) \sum_{|i| \ge 2} a_{i,i} - 2\sum_{|i| = 1,|j| \ge
2}a_{i,j} - 2\sum_{|i| \ge 2,|j| = 1}a_{i,j}.
\ena

Next, the equal contributions from $B_{I,J}={\bf 1}(\pi(I)=J,|I|
\ge 3)$ and $B_{J,I}$ sum to
\bea
\label{B3terms}
6 \sum_{|i| \ge 3} a_{i,\pi(i)} - 4 \sum_{|i| \ge 3}
a_{\pi^{-1}(i),\pi(i)} - 2 \sum_{|i| \ge 3} a_{\pi(i),i}.
\ena

On $F={\bf 1}(|I| \ge 2,|J| \ge 2,I \not = J, \pi(I)\not = J,
\pi(J) \not =I)$, the contribution from $a_{\pi^{-1}(I),I}$ is
\bea
\label{4-indicator}
\sum_{|i|,|j| \ge 2}
a_{\pi^{-1}(i),i}{\bf 1}(i \not = j, \pi(i)\not = j, \pi(j) \not
=i).
\ena
Let $i \cong j$ denote the fact that $i$ and $j$ are elements of
the same cycle. When $i \cong j$ and $\{i,j,\pi(i),\pi(j)\}$ are
distinct, we have $|i| \ge 4$ and $|i|-3$ possible choices for $j
\cong i$ that satisfy the conditions in the indicator in
(\ref{4-indicator}). Hence, the case $i \cong j$ contributes
\beas
\sum_{|i| \ge 4} a_{\pi^{-1}(i),i}\sum_{j \cong i}{\bf 1}(i \not =
j, \pi(i)\not = j, \pi(j) \not =i) = \sum_{|i| \ge 4}
a_{\pi^{-1}(i),i}(|i|-3)
 =\sum_{|i| \ge 3} (|i|-3) a_{i,\pi(i)}.
\enas
When $i \not \cong j$ the conditions in the indicator function
(\ref{4-indicator}) are satisfied if and only if $|i| \ge 2, |j|
\ge 2$. For $|i| \ge 2$ there are $n-|i|-c_1(\pi)$ choices for
$j$, so the case $i \not \cong j$ contributes
\beas
\sum_{|i| \ge 2} a_{\pi^{-1}(i),i}\sum_{j \not \cong i, |j| \ge
2}1 = \sum_{|i| \ge 2} (n-|i|-c_1(\pi)) a_{i,\pi(i)}.
\enas
The next three terms on $F$ give the same as the first, so in
total we have
\bea
\label{14-1terms}
4(n-2-c_1(\pi))\sum_{|i|=2}a_{i,\pi(i)}+ 4(n-3-c_1(\pi))\sum_{|i|
\ge 3}a_{i,\pi(i)}.
\ena

Decomposing the contribution from the fifth term, according to
whether $i \cong j$ or $i \not \cong j$, gives
\bea
\nonumber
&& -\sum_{|i|, |j| \ge 2} a_{\pi^{-1}(i),j}{\bf 1}(i
\not = j, \pi(i) \not = j, \pi(j) \not =i)\\
\nonumber && -\sum_{|i| \ge 4} \sum_{j \cong i} a_{\pi^{-1}(i),j}
{\bf 1}(i \not = j, \pi(i) \not = j, \pi(j) \not =i)- \sum_{|i|,|j| \ge 2} \sum_{j \not \cong i}
a_{\pi^{-1}(i),j}\\
&=& \nonumber -\sum_{|i| \ge 4} \sum_{j \cong i} a_{\pi^{-1}(i),j}
+ \sum_{|i| \ge 4}\left( a_{\pi^{-1}(i),i}+a_{\pi^{-1}(i),\pi(i)}
+ a_{\pi^{-1}(i),\pi^{-1}(i)} \right) - \sum_{|i|,|j| \ge 2} \sum_{j \not \cong i} a_{i,j}\\
\label{hope} &=& - \sum_{|i| \ge 4} \sum_{j \cong i} a_{i,j} +
\sum_{|i| \ge 4}\left( a_{i,\pi(i)}+a_{\pi^{-1}(i),\pi(i)} +
a_{i,i} \right)- \sum_{|i|,|j| \ge 2} \sum_{j \not \cong i}
a_{i,j}.
\ena

To simplify (\ref{hope}), let $a \wedge b = \min(a,b)$ and
consider the decomposition
\bea
\label{sumaij=0-miracle}
\sum_{i,j=1}^n a_{i,j}=\sum_{|i| \ge
4}\sum_{j \cong i}a_{i,j} + \sum_{|i| \le 3}\sum_{j \cong
i}a_{i,j}+ \sum_{|i|,|j| \ge 2} \sum_{j \not \cong i}a_{i,j}+
\sum_{|i| \wedge |j| = 1} \sum_{j \not \cong i}a_{i,j}.
\ena
Since $\sum_{i,j}a_{ij}=0$, we may replace the sum of the first
and last terms in (\ref{hope}) by the sum of the second and fourth
terms on the right hand side of (\ref{sumaij=0-miracle}),
respectively, resulting in
\beas
&& \sum_{|i| \le 3}\sum_{j \cong i}a_{i,j} + \sum_{|i| \wedge |j|
= 1} \sum_{j \not \cong i}a_{i,j} + \sum_{|i| \ge 4}\left(
a_{i,\pi(i)}+a_{\pi^{-1}(i),\pi(i)} +
a_{i,i}\right)\\
&=& \sum_{|i| \le 2}\sum_{j \cong i}a_{i,j} + \sum_{|i| \wedge |j|
= 1} \sum_{j \not \cong i}a_{i,j} + \sum_{|i| \ge 3}\left(
a_{i,\pi(i)}+a_{\pi^{-1}(i),\pi(i)} + a_{i,i}\right),
\enas
where we have used the fact that $\pi^2(i)=\pi^{-1}(i)$ when
$|i|=3$. Similarly shifting the $|i|=2$ term we obtain
\beas
&& \sum_{|i| = 1}a_{i,i} + \sum_{|i| \wedge |j| = 1} \sum_{j \not
\cong i}a_{i,j} + \sum_{|i| \ge 2}\left( a_{i,\pi(i)} +
a_{i,i}\right) + \sum_{|i| \ge 3}a_{\pi^{-1}(i),\pi(i)}\\
&=& \sum_{|i| \wedge |j| = 1} \sum_{j \not \cong i}a_{i,j} +
\sum_{|i| \ge 2} a_{i,\pi(i)} + \sum_{|i| \ge 1} a_{i,i}  +
\sum_{|i| \ge 3}a_{\pi^{-1}(i),\pi(i)}.
\enas
Combining this with the next three terms of $F$, each of which
yields the same contribution, gives
\bea
\label{C-terms} 4 \sum_{|i| \ge 2} a_{i,\pi(i)}+ 4 \sum_{|i| \ge
3} a_{\pi^{-1}(i),\pi(i)} + 4 \sum_{|i| \ge 1} a_{i,i} +
4\sum_{|i|\wedge |j| =1} \sum_{j \not \cong i}a_{i,j}.
\ena

Combining (\ref{C-terms}) with the contribution (\ref{14-1terms})
of the first four terms in $F$, the $A_{I,J}$ and $A_{J,I}$ terms
in (\ref{A3terms}) and the $B_{I,J}$ and $B_{J,I}$ terms
(\ref{B3terms}), yields $n(n-1)\E (Y'-Y''|\pi')$; after cancelling
the terms involving $a_{\pi^{-1}(i),\pi(i)}$ in (\ref{B3terms})
and (\ref{C-terms}) and grouping like terms, we obtain
\bea
\label{requires-symmetry}
&&4(n-1)\sum_{|i|=2}a_{i,\pi(i)}+(4n-2)\sum_{|i|\ge 3
}a_{i,\pi(i)} -2\sum_{|i|\ge 3}a_{\pi(i),i}\\
\label{requires-aii=0}
&+& 2(n-c_1(\pi)+2)\sum_{|i|=1}a_{i,i}-2(c_1(\pi)-2)\sum_{|i| \ge 2}a_{i,i}\\
\label{requires-no-1-cycles} &+&4 \sum_{|i| \wedge |j|=1, j \not
\cong i}a_{i,j} -2\sum_{|i|=1,|j|\ge 2}a_{i,j} -2\sum_{|i| \ge
2,|j|=1}a_{i,j}.
\ena
The assumption that $a_{i,i}=0$ causes the contribution from
(\ref{requires-aii=0}) to vanish, the assumption that there are no
1-cycles causes the contribution from (\ref{requires-no-1-cycles})
to vanish, and the assumption that $a_{i,j}$ is symmetric causes
the combination of the second and third terms in
(\ref{requires-symmetry}) to yield $\E
(Y'-Y''|\pi')=(4/n)\sum_{i=1}^n a_{i,\pi'(i)}=(4/n)Y'$. Hence, the
linearity condition (\ref{linear-lambda-F}) is satisfied.

Since $\pi''=\tau_{IJ} \pi \tau_{IJ}$, the terms that multiply the
indicator functions in the difference $Y'-Y''$ in (\ref{A}) depend
only on values in a subset of
$\{\pi^{-1}(I),I,\pi(I),\pi^{-1}(J),J, \pi(J)\}$ determined by the
event indicated; for example, on $B_{I,J}$ the difference only
depends on $\{ \pi^{-1}(I),I,J,\pi(J)\}$. For each event we
tabulate such values in a vector ${\bf i}$. Likewise, with
$\pi^\dagger$ and $\pi^\ddagger$ constructed according to
$\pi^\ddagger=\tau_{I^\dagger J^\dagger} \pi^\dagger
\tau_{I^\dagger J^\dagger}$, the difference $Y^\dagger-Y^\ddagger$
depends only on a subset of $\{P^\dagger,I^\dagger,K^\dagger
,Q^\dagger,J^\dagger, L^\dagger\}$, the corresponding values in
the $\pi^\dagger$ cycle, which we will tabulate in a vector ${\bf
i}^\dagger$. Since $Y'-Y''$ in (\ref{A}) is a sum of terms
multiplied by indicator functions of disjoint events, $(Y'-Y'')^2$
is a sum of those terms squared, multiplied by the same indicator
functions. Hence to generate $(\pi^\dagger,\pi^\ddagger)$ such
that $(Y^\dagger,Y^\ddagger)$ has a distribution proportional to
$(y'-y'')^2dF(y',y'')$, on each event we generate the elements of
${\bf i}^\dagger$ with square weighted probability appropriate to
the set indicated. Once the values in ${\bf i}^\dagger$ are
chosen, in order for $\pi^\dagger$ to have the conditional
distribution of $\pi$ given these values, the remaining values of
$\pi^\dagger$ are obtained by interchanging ${\bf i}$ with ${\bf
i}^\dagger$ in the cycle structure of $\pi$. That is, in each case
we specify $\pi^\dagger$ in terms of $\pi$ by
\bea
\label{def-pi-dagger-transposition} \pi^\dagger = \tau_{{\bf
i},{\bf i}^\dagger} \pi \tau_{{\bf i},{\bf i}^\dagger} \quad
\mbox{where} \quad \tau_{{\bf i},{\bf
i}^\dagger}=\prod_{k=1}^{\kappa} \tau_{i_k,i_k^\dagger},
\ena
and ${\bf i}=(i_1, \ldots, i_\kappa)$ and ${\bf
i}^\dagger=(i_1^\dagger,\ldots,i_\kappa^\dagger)$ are vectors of
disjoint indices, of some length $\kappa$.

For $\rho \in {\cal S}_\kappa$ and ${\bf l}=(l_1,\ldots,l_\kappa)$
any $\kappa$-dimensional vector of indices, let $\rho({\bf
l})=\{\rho(l_k):k=1,\ldots,\kappa\}$, and let $\iota$ be the
identity permutation. Since the values of
$\tau_{ii^\dagger}\pi\tau_{ii^\dagger}$ may differ from those of
$\pi$ only at $i,i^\dagger,\pi^{-1}(i)$ and $\pi^{-1}(i^\dagger)$,
(\ref{new-def-T}) will hold for the variables given by
(\ref{def-T'}), with
\beas
{\cal I}= \iota({\bf i}) \cup \iota({\bf i}^\dagger) \cup
\pi^{-1}({\bf i}) \cup \pi^{-1}({\bf i}^\dagger).
\enas

The construction in each case proceeds as follows. Since 1-cycles
are excluded, $A_{I,J}$ and $A_{J,I}$ are null. On $B_{I,J}$,
where $I,J$ and $\pi(J)$ are three distinct, consecutive values of
$\pi$, if $|I|=3$ then the symmetry of $a_{i,j}$ gives $Y''=Y'$,
an event on which the distribution of $(Y^\dagger,Y^\ddagger)$,
proportional to $(Y''-Y')^2$, puts mass zero. Otherwise, $|I| \ge
4$ and $Y'-Y''$ depends only on ${\bf
i}=(\pi^{-1}(I),I,J,\pi(J))$, and we choose ${\bf
i}^\dagger=(P^\dagger,I^\dagger,J^\dagger, L^\dagger)$, the
corresponding values for $\pi^\dagger$, according to the
distribution
\beas
(P^\dagger,I^\dagger,J^\dagger,L^\dagger) &\sim&
[(a_{p,i}+a_{j,l})-(a_{p,j}+a_{i,l})]^2{\bf 1}(\mbox{$p,i,j$, and
$l$ are distinct}),
\enas
noting that $a_{i,j}$ cancels with $a_{j,i}$ by symmetry. Now set
$\pi^\dagger$ as specified in (\ref{def-pi-dagger-transposition}).
In this case ${\cal I}$ has size at most thirteen. Reversing the
roles of $I$ and $J$ gives the construction on $B_{J,I}$.

Next consider $F$, where $I,\pi(I),J$ and $\pi(J)$ are distinct.
If $|I|=|J|=2$ then take
\beas
(I^\dagger,K^\dagger,J^\dagger,L^\dagger) \sim [(a_{i,k}+a_{j,l})-
(a_{i,l}+a_{j,k})]^2 {\bf 1}(\mbox{$i,k,j$ and $l$ are distinct}),
\enas
and set $\pi^\dagger$ as specified in
(\ref{def-pi-dagger-transposition}), with ${\bf
i}=(I,\pi(I),J,\pi(J))$ and ${\bf
i}^\dagger=(I^\dagger,K^\dagger,J^\dagger, L^\dagger)$, and with
the size of ${\cal I}$ at most twelve. For $|I| \ge 3$ and
$|J|=2$, take
\beas
(P^\dagger,I^\dagger, K^\dagger, J^\dagger, L^\dagger) \sim
[(a_{p,i}+ a_{i,k}+2a_{j,l})- (a_{p,j}+a_{j,k}+2a_{i,l})]^2 {\bf
1}(\mbox{$p,i,k,j,$ and $l$ are distinct}),
\enas
and set $\pi^\dagger$ as specified in
(\ref{def-pi-dagger-transposition}), with ${\bf
i}=(\pi^{-1}(I),I,\pi(I),J,\pi(J))$ and ${\bf
i}^\dagger=(P^\dagger,I^\dagger,K^\dagger,J^\dagger,L^\dagger)$,
and with the size of ${\cal I}$ at most sixteen. Reversing the
roles of $I$ and $J$ gives the case in which $|J|=2$ but $|I|\ge
3$. For $|I| \ge 3, |J| \ge 3$, take
\beas (P^\dagger,I^\dagger,K^\dagger,Q^\dagger,J^\dagger,L^\dagger)
&\sim& [(a_{p,i}+a_{i,k}+a_{q,j}+ a_{j,l})-
(a_{p,j}+a_{j,k}+a_{q,i}+a_{i,l})]^2 \\
&\times&{\bf 1}(\mbox{$p,i,k,q,j,$ and $l$ distinct}),
\enas and set $\pi^\dagger$ as specified in
(\ref{def-pi-dagger-transposition}), with ${\bf
i}=(\pi^{-1}(I),I,\pi(I),\pi^{-1}(J),J,\pi(J))$ and \\${\bf
i}^\dagger=(P^\dagger,I^\dagger,K^\dagger,Q^\dagger,J^\dagger,L^\dagger)$.
In this case, the size of ${\cal I}$ is at most twenty and, by
(\ref{Y*Y2B-bound}), $|Y^*-Y| \le 40C$ in all cases. $\qed$

\section{Size Biasing: Permutations and Patterns}
\label{size-applications}

In this section we derive corollaries of Theorem
\ref{size-delta-bound-theorem} to obtain Berry Esseen bounds for
the number of occurrences of fixed, relatively ordered
sub-sequences, such as rising sequences, in a random permutation,
and of color patterns, local maxima, and sub-graphs in finite
graphs.

Following [\ref{Goldstein-Rinott}], given a finite collection
${\bf X}=\{X_\alpha, \alpha \in {\cal A}\}$ of non-negative random
variables with index set ${\cal A}$, for $\alpha \in {\cal A}$ we
say the collection ${\bf X}^\alpha=\{X_\beta^\alpha, \beta \in
{\cal A}\}$ has the ${\bf X}$-size-biased distribution in
direction $\alpha$ if
\bea
\label{size-bias-alpha} \E X_\alpha f({\bf X}) = \E X_\alpha \E
f({\bf X}^\alpha)
\ena
for all functions $f$ on ${\bf X}$ for which these expectations
exist. For the given ${\bf X}$, the collection ${\bf X}^\alpha$
exists for any $\alpha \in {\cal A}$ and has distribution
$dP^\alpha({\bf x})=x_\alpha dP({\bf x})/\E X_\alpha$, where
$dP({\bf x})$ is the distribution of ${\bf X}$. Specializing
(\ref{size-bias-alpha}) to the coordinate function $f({\bf
X})=g(X_\alpha)$, we see that $X_\alpha^\alpha$ has the
$X_\alpha$-size-biased distribution $X_\alpha^s$, defined in
(\ref{size-bias}).

\begin{corollary}
\label{local-size-bias} Let ${\bf X}=\{X_\alpha,\alpha \in {\cal
A}\}$ be a finite collection of random variables with values in
$[0,M]$ and let $Y=\sum_{\alpha \in {\cal A}}X_\alpha$. Assume,
for each $\alpha \in {\cal A}$, there exists a dependency
neighborhood ${\cal B}_\alpha \subset {\cal A}$ such that
\bea
\label{X-alpha-indep-of-Balpha} X_\alpha \quad \mbox{and} \quad
\{X_\beta: \beta \not \in {\cal B}_\alpha \} \quad \mbox{are
independent.}
\ena
Furthermore, let $p_\alpha=\E  X_\alpha/\sum_{\beta \in {\cal A}}
\E  X_\beta$ and $\max_\alpha |{\cal B}_\alpha| = b$. For each
$\alpha \in {\cal A}$, let $({\bf X},{\bf X}^\alpha)$ be a
coupling of ${\bf X}$ to an ${\bf X}^\alpha$ with the ${\bf
X}$-size-biased distribution in direction $\alpha$, and let ${\cal
D} \subset {\cal A} \times {\cal A}$ and ${\cal F} \supset
\sigma\{Y\}$ be such that if $(\alpha_1,\alpha_2) \not \in {\cal
D}$ then
\bea
\label{cov=0} \mbox{Cov}(\E (X_{\beta_1}^{\alpha_1}-
X_{\beta_1}|{\cal F}),\E (X_{\beta_2}^{\alpha_2}-X_{\beta_2}|{\cal
F}))=0  \quad \mbox{for all $(\beta_1,\beta_2) \in {\cal
B}_{\alpha_1} \times  {\cal B}_{\alpha_2}$.}
\ena
Then Theorem \ref{size-delta-bound-theorem} may be applied with
\bea
\label{for-graph-bounds} B=bM \quad \mbox{and} \quad \Delta \le
M\sqrt{\sum_{(\alpha_1,\alpha_2) \in {\cal
D}}p_{\alpha_1}p_{\alpha_2}|{\cal B}_{\alpha_1}||{\cal
B}_{\alpha_2}|} \le (\max_\alpha p_\alpha) bM \sqrt{|{\cal D|}}.
\ena
\end{corollary}

\noindent {\bf Proof:} Assuming, without loss of generality, that
$\E X_\alpha >0$ for each $\alpha \in {\cal A}$, the factorization
\beas
P^\alpha({\bf X}\in d{\bf x}) =
\left( \frac{ x_\alpha P(X_\alpha \in dx_\alpha)}{ \E  X_\alpha}
\right) P({\bf X} \in d{\bf x}\,| \,X_\alpha = x_\alpha)
\enas
shows that we can construct ${\bf X}^\alpha$ by first choosing
$X_\alpha^\alpha$ from the $X_\alpha$-size-bias distribution, and
then choosing the remaining variables from the conditional
distribution of ${\bf X}$, given the chosen value of
$X_\alpha^\alpha$. Note that $X_\beta^\alpha \in [0,M]$ for all
$\alpha,\beta$ and, by (\ref{X-alpha-indep-of-Balpha}), that we
may take $X_\beta^\alpha=X_\beta$ for $\beta \not \in {\cal
B}_\alpha$. By Lemma 2.1 of [\ref{Goldstein-Rinott}],
$Y^s=\sum_{\beta \in {\cal A}}X_\beta^I$ has the $Y$-size-biased
distribution, where the random index $I$ has distribution
$P(I=\alpha)=p_\alpha$, and is independent of both $({\bf X},{\bf
X}^\alpha)$ and $\cal{F}$. Hence
\bea
\label{Ys-Y} Y^s-Y=\sum_{\beta \in {\cal B}_I}(X_\beta^I-X_\beta),
\quad \mbox{and therefore,} \quad |Y^s-Y| \le b M.
\ena

Since $\sigma\{Y\} \subset {\cal F}$, $ \Delta^2 = \mbox{Var}(\E
(Y^s-Y|Y)) \le \mbox{Var}(\E  (Y^s-Y| {\cal F})). $ Taking
conditional expectation with respect to ${\cal F}$ in (\ref{Ys-Y})
yields,
$$
\E  \left(Y^s-Y| {\cal F}\right) = \sum_{\alpha \in {\cal A}}
p_\alpha \sum_{\beta \in {\cal B}_\alpha}\E (X_\beta^\alpha -
X_\beta|{\cal F})
$$
and, therefore,
\beas
\mbox{Var}(\E  \left(Y^s-Y| {\cal F}\right)) &=& \E
\sum_{\stackrel{(\alpha_1,\alpha_2) \in {\cal A} \times {\cal A}}
{(\beta_1,\beta_2) \in {\cal B}_{\alpha_1} \times {\cal
B}_{\alpha_2}}} p_{\alpha_1} p_{\alpha_2} \mbox{Cov}(\E
(X_{\beta_1}^{\alpha_1}- X_{\beta_1}|{\cal F}),\E
(X_{\beta_2}^{\alpha_2}-X_{\beta_2}|{\cal F})).
\enas
Using (\ref{cov=0}), we may replace the sum over
$(\alpha_1,\alpha_2) \in {\cal A} \times {\cal A}$ by the sum over
$(\alpha_1,\alpha_2) \in {\cal D}$, and subsequent application of
the Cauchy Schwarz inequality yields the bound
(\ref{for-graph-bounds}) for $\Delta$. $\qed$

If, in some asymptotic regime, the $X_\alpha$ are comparable in
expectation in such a way that $p_\alpha \sim |{\cal A}|^{-1}$; if
$\mu$ and $\sigma^2$ grow like $|{\cal A}|$; if $b$ remains
bounded; and if $|{\cal D}|$ is of order $|{\cal A}|$, then, in
Theorem \ref{size-delta-bound-theorem}, $A$ and $\Delta$ and,
therefore, $\delta$ are of order $1/\sigma$.

\begin{corollary}
\label{apple-graph-corollary} Let ${\cal G}$ be an index set, let
$\{C_g, g \in {\cal G}\}$ be a collection of independent random
elements taking values in an arbitrary set ${\cal C}$, let
$\{{\cal G}_\alpha, \alpha \in {\cal A}\}$ be a finite collection
of subsets of ${\cal G}$, and, for $\alpha \in {\cal A}$, let
\beas
X_\alpha=X_\alpha(C_g: g \in {\cal G}_\alpha)
\enas
be a function of the variables $\{C_g, g \in {\cal G}_\alpha\}$,
taking values in $[0,M]$. Then Theorem
\ref{size-delta-bound-theorem} may be applied to $Y=\sum_\alpha
X_\alpha$ with $B$ and $\Delta$ as in (\ref{for-graph-bounds}),
taking $p_\alpha=\E X_\alpha/ \sum_\beta \E X_\beta$,
\bea
\label{B-independent} {\cal B}_\alpha = \{\beta \in {\cal A}:
{\cal G}_\beta \cap {\cal G}_\alpha \not = \emptyset\} \quad
\mbox{for $\alpha \in {\cal A},$}
\ena
and any ${\cal D}$ for which
\bea \label{D-independent}
{\cal D} \supset \{(\alpha_1,\alpha_2): \mbox{there exists
$(\beta_1,\beta_2) \in {\cal B}_{\alpha_1} \times {\cal
B}_{\alpha_2}$ with ${\cal G}_{\beta_1} \cap {\cal G}_{\beta_2}
\not = \emptyset$} \}.
\ena
\end{corollary}

\noindent {\bf Proof:} Since $X_\alpha$ and $X_\beta$ are
functions of disjoint sets of independent variables when ${\cal
G}_\alpha \cap {\cal G}_\beta = \emptyset$,
(\ref{X-alpha-indep-of-Balpha}) holds with the dependency
neighborhoods given by (\ref{B-independent}). Now, for each
$\alpha \in {\cal A}$, consider the following $({\bf X},{\bf
X}^\alpha)$ coupling. Let $\{C_g^{(\alpha)}, g \in {\cal
G}_\alpha\}$ be independent of $\{C_g, g \in {\cal G} \}$ and have
distribution
$$
dP^{(\alpha)}(c_g,g \in {\cal G}_\alpha) = \frac{X_\alpha(c_g,g
\in {\cal G}_\alpha)}{\E  X_\alpha(c_g,g \in {\cal
G}_\alpha)}dP(c_g,g \in {\cal G}_\alpha).
$$
Then, by direct verification of (\ref{size-bias-alpha}), the
collection
$$
X_\beta^\alpha = X_\beta(C_g, g \in {\cal G}_\beta \cap {\cal
G}_\alpha^c, \,\,C_g^{(\alpha)}, g \in {\cal G}_\beta \cap {\cal
G}_\alpha), \quad \beta \in {\cal A}
$$
has the ${\bf X}^\alpha$ distribution. Taking ${\cal F} = \{C_g: g
\in {\cal G}\}$, we have $\E (X_\beta^\alpha|{\cal F})=\E
(X_\beta^\alpha|C_g,g \in {\cal G}_\beta)$ and, since $\E
(X_\beta|{\cal F})= X_\beta$, the conditional expectation $\E
(X_\beta^\alpha-X_\beta|{\cal F})$ is a function of $\{C_g,g \in
{\cal G }_\beta \}$ only. In particular, if $(\alpha_1,\alpha_2)
\not \in {\cal D}$ then, for all $\beta_1 \in {\cal B}_{\alpha_1}$
and $\beta_2 \in {\cal B}_{\alpha_2}$ we have ${\cal G}_{\beta_1}
\cap {\cal G}_{\beta_2} = \emptyset$ and, consequently, $\E
(X_{\beta_1}^{\alpha_1}-X_{\beta_1}|{\cal F})$ and $\E
(X_{\beta_2}^{\alpha_2}-X_{\beta_2}|{\cal F})$ are independent,
yielding (\ref{cov=0}), and all conditions of Corollary
\ref{local-size-bias} hold. $\qed$

With the exception of Example \ref{perm}, in the remainder of this
section we consider graphs ${\cal G}=({\cal V},{\cal E})$ having
random elements $\{C_g\}_{g \in {{\cal V} \cup {\cal E}}}$
assigned to their vertices and edges, and applications of
Corollary \ref{apple-graph-corollary} to the sum $Y=\sum_{\alpha
\in {\cal A}} X_\alpha$ of bounded functions
$X_\alpha=X_\alpha(C_g, g \in {\cal V}_\alpha \cup {\cal
E}_\alpha)$, where ${\cal G}_\alpha=({\cal V}_\alpha, {\cal
E}_\alpha), \alpha \in {\cal A}$ is a given finite family of
subgraphs of ${\cal G}$; we abuse notation slightly in that a
graph ${\cal G}$ is replaced by ${\cal V} \cup {\cal E}$ when used
as an index set for the underlying variables $C_g$. When
$\{C_g\}_{g \in {\cal G}}$ are independent, Corollary
\ref{apple-graph-corollary} applies and, in (\ref{B-independent})
and (\ref{D-independent}), the intersection of the two graphs
$({\cal V}_1,{\cal E}_1)$ and $({\cal V}_2,{\cal E}_2)$ is the
graph $({\cal V}_1 \cap {\cal V}_2, {\cal E}_1 \cap {\cal E}_2)$.

Furthermore, if ${\cal A} \subset {\cal V}$ and there is a
distance $d(\alpha,\beta)$ defined on ${\cal A}$, then letting
\bea
\label{def-rho}
\rho = \inf\{\varrho: {\cal V}_\alpha \cap {\cal
V}_\beta = \emptyset\, \mbox{for all $\alpha,\beta \in {\cal A}$
with $d(\alpha,\beta)
> \varrho$} \},
\ena
we may use
\bea
\label{B-D-rho-regular} {\cal B}_\alpha = \{\beta: d(\alpha,\beta)
\le \rho \} \quad \mbox{and} \quad {\cal D}=\{(\alpha_1,\alpha_2):
d(\alpha_1,\alpha_2) \le 3 \rho \}
\ena
in (\ref{B-independent}) and (\ref{D-independent}), respectively,
since rearranging $d(\alpha_1,\alpha_2) \le
d(\alpha_1,\beta_1)+d(\beta_1,\beta_2)+d(\beta_2,\alpha_2)$ gives,
\beas
d(\beta_1,\beta_2) \ge
d(\alpha_1,\alpha_2)-(d(\alpha_1,\beta_1)+d(\alpha_2,\beta_2)) \ge
d(\alpha_1,\alpha_2)-2\rho > \rho,
\enas
for $(\alpha_1,\alpha_2) \not \in {\cal D}$ and $(\beta_1,\beta_2)
\in {\cal B}_{\alpha_1} \times {\cal B}_{\alpha_2}$.

For $v \in {\cal V}$ and $r \ge 0$ let ${\cal G}_{v,r}$ be the
restriction of ${\cal G}$ to the vertices at most a distance $r$
from $v$; that is ${\cal G}_{v,r}$ has vertex set ${\cal
V}_{v,r}=\{w \in {\cal V}: d(v,w) \le r \}$ and edge set ${\cal
E}_{v,r}= \{\{w,u\} \in {\cal E}: w,u \in {\cal V}_{v,r}\}$. We
say that a graph ${\cal G}$ is {\em distance $r$-regular} if
${\cal G}_{v,r}$ is isomorphic to some graph $({\cal V}_r,{\cal
E}_r)$ for all $v$. For example, a graph of constant degree is
distance 1-regular. This notion of distance $r$-regular is related
to, but not the same as, the notion of a distance-regular graph as
given in [\ref{Biggs}] and [\ref{Brouwer}]. For a distance
$r$-regular graph let
\bea
\label{def-Vr}
V(r)=|{\cal V}_r|.
\ena
Corollary \ref{distance-regular}, below, follows from Corollary
\ref{apple-graph-corollary} as a consequence of the remarks above,
and by noting that the given assumptions imply that $|{\cal
D}|=|{\cal A}|V(3\rho)$ and that $\E X_\alpha$ is constant,
yielding $p_\alpha=1/|{\cal A}|$.

\begin{corollary}
\label{distance-regular} Let ${\cal G}$ be a graph with a finite
family of isomorphic subgraphs $\{{\cal G}_\alpha,\alpha \in {\cal
A}\}, {\cal A} \subset {\cal V}$, let $d(\cdot,\cdot)$ be a
distance on ${\cal A}$, and define $\rho$ as in (\ref{def-rho}).
For each $\alpha \in {\cal A}$, let $X_\alpha$ be given by
\bea
\label{same-X} X_\alpha=X(C_g, g \in {\cal G}_\alpha)
\ena
for a fixed function $X$ taking values in $[0,M]$, and let the
elements of $\{C_g\}_{g \in {\cal G}}$ be independent, with
$\{C_g: g \in {\cal G}_\alpha\}$ identically distributed. If
${\cal G}$ is a distance-$3 \rho$-regular graph, then Theorem
\ref{size-delta-bound-theorem} may be applied to $Y=\sum_{\alpha
\in {\cal A}} X_\alpha$ with $V(r)$ as given in (\ref{def-Vr}) and
\bea
\label{cor-dist-reg-A-Delta} B=V(\rho)M, \quad \Delta \le M|{\cal
A}|^{-1/2}V(\rho)\sqrt{V(3\rho)}.
\ena
\end{corollary}

Natural families of examples in $\mathbb{R}^p$ can be generated
using the vertex set ${\cal V}=\{1,\ldots,n\}^p$ with
componentwise addition modulo $n$, and $d(\alpha,\beta)$ given by
e.g. the $L^1$ distance $||\alpha-\beta||$.

\begin{example}
\label{m-dependent-example} ({\em Sliding $m$-window.}) For $n \ge
m \ge 1$, let ${\cal A}={\cal V}=\{1,\ldots,n\}$ considered modulo
$n$, $\{C_g:g \in {\cal G}\}$ i.i.d. real valued random variables,
and for each $\alpha \in {\cal A}$
\bea
\label{calG-oned} {\cal G}_\alpha= ({\cal V}_\alpha,{\cal
E}_\alpha), \quad \mbox{where} \quad {\cal V}_\alpha=\{v \in {\cal
V}: \alpha \le v \le \alpha + m -1 \} \quad \mbox{and} \quad {\cal
E}_\alpha=\emptyset.
\ena
Then for $X:\mathbb{R}^m \rightarrow [0,1]$, Corollary
\ref{distance-regular} may be applied to the sum $Y=\sum_{\alpha
\in {\cal A}}X_\alpha$ of the $m$-dependent sequence $X_\alpha =
X(C_\alpha, \ldots, C_{\alpha+m-1})$, formed by applying the
function $X$ to the variables in the `$m$-window' ${\cal
V}_\alpha$. In this example, taking
$d(\alpha,\beta)=|\alpha-\beta|$ gives $\rho=m-1$ and $V(r)=2r+1$.
Hence, from (\ref{cor-dist-reg-A-Delta}), $B=(2m-1)$ and $\Delta
\le n^{-1/2}(2m-1)(6m-5)^{1/2}$.
\end{example}

In Example \ref{perm} the underlying variables are not
independent, and Corollaries \ref{apple-graph-corollary} and
\ref{distance-regular} cannot be directly applied.
\begin{example}
\label{perm} ({\em Relatively ordered sub-sequences of a random
permutation.}) For $n \ge m \ge 1$, let $\pi$ be a uniform random
permutation of the integers ${\cal V}=\{1,\ldots,n\}$, taken
modulo $n$. For a permutation $\tau$ on $\{1,\ldots,m \}$, let
${\cal G}_\alpha$ and ${\cal V}_\alpha$ be as specified in
(\ref{calG-oned}), and let $X_\alpha$ the indicator function
requiring that the pattern $\tau$ appears on ${\cal V}_\alpha$;
that is, that the values $\{\pi(v)\}_{v \in {\cal V}_\alpha}$ and
$\{\tau(v)\}_{v \in {\cal V}_1}$ are in the same relative order.
Equivalently, the pattern $\tau$ appears on ${\cal V}_\alpha$ if
and only if $\pi(\tau^{-1}(v)+\alpha-1), v  \in {\cal V}_1$ is an
increasing sequence, and we write
$$
X_\alpha(\pi(v), v \in {\cal G}_\alpha) = {\bf
1}(\pi(\tau^{-1}(1)+\alpha -1 ) < \cdots <
\pi(\tau^{-1}(m)+\alpha-1)).
$$
With ${\cal A}={\cal V}$, the sum $Y=\sum_{\alpha \in {\cal A}}
X_\alpha$ counts the number of $m$-element-long segments of $\pi$
that have the same relative order as $\tau$.

For $\alpha \in {\cal A}$, we generate ${\bf
X}^\alpha=\{X_\beta^\alpha, \beta \in {\cal A}\}$ by reordering
the values of $\pi(\gamma)$ for $\gamma \in {\cal V}_\alpha$, to
be in the same relative order as $\tau$, and let $X_\beta^\alpha$
be the indicator requiring $\tau$ to appear at position $\beta$ in
the reordered permutation. Letting ${\cal F}=\sigma\{\pi\}$, we
have $\E (X_\beta^\alpha|{\cal F})$ and $X_\beta$ depend only on
the relative order of $\{\pi(\gamma): -(m-1) \le \gamma - \beta
\le 2(m-1)\}$. Since the relative order of the non-overlapping
segments of the values of $\pi$ are independent,
(\ref{X-alpha-indep-of-Balpha}) and (\ref{cov=0}) hold when ${\cal
B}_\alpha$ and ${\cal D}$ are as in (\ref{B-D-rho-regular}), for
$d(\alpha,\beta)=|\alpha-\beta|$ and $\rho=m-1$; hence, Theorem
\ref{size-delta-bound-theorem} may be applied with the same value
for $B$ and bound on $\Delta$ as in Example
\ref{m-dependent-example}.

When $\tau=\iota_m$, the identity permutation of length $m$, we
say that $\pi$ has a rising sequence of length $m$ at position
$\alpha$ if $X_\alpha=1$. Rising sequences were studied in
[\ref{Dia}] in connection with card tricks and card shuffling. Due
to the regular-self-overlap property of rising sequences, namely
that a non-empty intersection of two rising sequences is again a
rising sequence, some improvement on the constant in the bound can
be obtained by a more careful consideration of the conditional
variance.
\end{example}

\begin{example}
({\em Coloring patterns and subgraph occurrences on a finite graph
$\,{\cal G}$}). For illustration, take ${\cal V}={\cal
A}=\{1,\ldots,n\}^p$, considered modulo $n$, let
$d(\alpha,\beta)=||\alpha-\beta||$ with $||\cdot||$ the sup norm,
let ${\cal E}=\{\{w,v\}: d(w,v)=1\}$, and, for each $\alpha \in
{\cal A}$, let ${\cal G}_\alpha=({\cal V}_\alpha,{\cal E}_\alpha)$
where
\beas
 {\cal V}_\alpha = \{\alpha + (e_1,\ldots,e_p): e_i
\in \{0,1\}\} \quad \mbox{and} \quad {\cal E}_\alpha=\{ \{v,w\}:
v,w \in {\cal V}_\alpha,\,\, d(w,v)=1\}.
\enas
Let ${\cal C}$ be a set (of e.g. colors) from which is formed a
given pattern $\{c_g: g \in {\cal G}_{\bf 0}\}$, let $\{C_g, g \in
{\cal G}\}$ be independent variables in ${\cal C}$ with $\{C_g: g
\in {\cal G}_\alpha \}_{\alpha \in {\cal A}}$ identically
distributed, and let
\bea
\label{X-color-subgraph} X(C_g,g \in {\cal G}_{\bf 0})=\prod_{g
\in {\cal G}_{\bf 0}}{\bf 1}(C_g = c_g),
\ena
and $X_\alpha$ given by (\ref{same-X}). Then $Y=\sum_{\alpha \in
{\cal A}} X_\alpha$ counts the number of times the pattern appears
in the subgraphs ${\cal G}_\alpha$. Corollary
\ref{distance-regular} may be applied with $M=1$, $\rho=1$ (by
(\ref{def-rho})), $V(r)=(2r+1)^p$, and (by
(\ref{cor-dist-reg-A-Delta})) $B=3^p$ and $\Delta \le
(63/n)^{p/2}$.

 Such multi-dimensional pattern occurrences are a
generalization of the well-studied case in which one-dimensional
sequences are scanned for pattern occurrences; see, for instance,
[\ref{naus-book}] and [\ref{naus-scan}] for scan and window
statistics, see [\ref{huang}] for applications of the normal
approximation in this context to molecular sequence data, and see
also [\ref{DW-1}] and [\ref{DW-2}], where higher-dimensional
extensions are considered.

Occurrences of subgraphs can be handled as a special case. For
example, with $({\cal V},{\cal E})$ the graph above, let $G$ be
the random subgraph with vertex set ${\cal V}$ and random edge set
$\{e \in {\cal E}: C_e=1\}$ where $\{C_e\}_{e \in {\cal E}}$ are
independent and identically distributed Bernoulli variables. Then
say, taking the product in (\ref{X-color-subgraph}) over edges $e
\in {\cal E}_0$ and setting $c_e=1$, the sum $Y=\sum_{\alpha \in
{\cal A}} X_\alpha$ counts the number of times that copies of
${\cal E}_0$ appear in the random graph $G$; the same bounds hold
as above.

The authors of [\ref{Barbour}] studied the related problem of
counting the number of small cliques that occur in the random
binomial graph, a case in which the dependence is not local; the
technique applied is the Chen-Stein method.
\end{example}

\begin{example} ({\em Local extremes.}) Let ${\cal G}_\alpha, \alpha \in
{\cal A}$, be a collection of subgraphs of ${\cal G}$ isomorphic
to ${\cal G}_0$, let $v \in {\cal V}_0$ be a distinguished vertex,
let $\{C_g, g \in {\cal V} \}$ be a collection of independent and
identically distributed random variables, and let $X_\alpha$ be
defined by (\ref{same-X}) with
$$
X(C_\beta, \beta \in {\cal V}_0) = {\bf 1}(C_v \ge C_\beta, \beta
\in {\cal V}_0).
$$
Then the sum $Y=\sum_{\alpha \in {\cal A}} X_\alpha$ counts the
number of times the vertex in ${\cal G}_\alpha$ which corresponds
under the isomorphism to the distinguished vertex $v \in {\cal
V}_0$, is a local maxima. Corollary \ref{distance-regular} holds
with $M=1$; the other quantities determining the bound begin
dependent on the structure of ${\cal G}$.

For example, consider the hypercube ${\cal V}=\{0,1\}^p$ and
${\cal E}=\{\{v,w\}:||v-w||=1\}$, where $||\cdot||$ is the Hamming
distance (see also [\ref{BRS}] and [\ref{BR}]). Take $v={\bf 0}$,
${\cal A}={\cal V}$, and, for each $\alpha \in {\cal A}$, let
${\cal V}_\alpha = \{\beta: ||\beta-\alpha|| \le 1 \}$ and ${\cal
E}_\alpha=\{\{v,w\}:v,w \in {\cal V}_\alpha, ||v-w||=1\}$.
Corollary \ref{distance-regular} applies with $\rho=2$ (by
(\ref{def-rho})), $V(r)=\sum_{j=0}^r {p \choose j}$, and (by
(\ref{cor-dist-reg-A-Delta}))
$$
B=1+p+{p \choose 2}\quad \mbox{and} \quad \Delta \le
2^{-p/2}\sum_{j=0}^2 {p \choose j} \sqrt{ \sum_{j=0}^6 {p \choose
j}}.
$$
\end{example}

\section{Proofs of Theorems 1.1 and 1.2}
\label{proofs} In this section, ${\cal H}$ denotes a class of
measurable functions satisfying properties (i),(ii), and (iii) (as
described in Section \ref{introduction}), and $h$ denotes an
element of ${\cal H}$. Recall that $\delta$ is given by
(\ref{def-delta}), let $\phi(t)$ denote the standard normal
density, and, for $t \in (0,1)$, define \beq \label{def-ht} h_t(x)
= \int h(x+ty) \phi(y) dy \quad \mbox{and} \quad \delta_t = \sup
\{ |\E h_t(W)-Nh_t|: h \in {\cal H} \}. \enq
\begin{lemma}
\label{bhat-smoothing} For a random variable $W$ on $\mathbb{R}$,
we have
\bea
\label{smooth}
\delta \leq 2.8 \delta_t + 4.7 at \quad \mbox{for
all $t \in (0,1)$,}
\ena
where $a$ is as in (\ref{aaa}). Furthermore, for all $A >0$ and
$\tilde{h}_\epsilon$ as in (\ref{H-closure}),
\bea
\label{2delta+a} \E  \left( \int \tilde{h}_{A +t|y|}(W) |\phi'(y)|
dy \right) \le 2\delta + a(A+t).
\ena
\end{lemma}

\noindent {\bf Proof:} Inequality (\ref{smooth}) is Lemma 4.1 of
[\ref{Rin-Rot-2}], following Lemma 2.11 of [\ref{Gotze}], which
stems from [\ref{bhat}]. As in [\ref{Rin-Rot-2}], adding and
subtracting to the left hand side of (\ref{2delta+a}) we have
\bea
\nonumber &&\E  \left( \int  (\tilde{h}_{A +t|y|}(W)- \tilde{h}_{A
+t|y|}(Z)) \,|\phi'(y)| dy + \int
\tilde{h}_{A +t|y|}(Z) |\phi'(y)|dy\right)\\
\label{apply-CS} &\le& \int  |\E  \tilde{h}_{A +t|y|}(W)- \E
\tilde{h}_{A +t|y|}(Z)| \,|\phi'(y)| dy + \int \E \tilde{h}_{A
+t|y|}(Z) |\phi'(y)|dy \\
\nonumber &\le& \left(2\delta + \int a(A+t|y|)|\phi'(y)|dy
\right)\le 2\delta + a(A+t),
\ena
where for the first term inside the parentheses in
(\ref{apply-CS}), we have used the facts that $h^{\pm}_{A +t|y|}
\in {\cal H}$ and $\int |\phi'(y)|dy \le 1$. For the second term
in the parentheses, we have used (\ref{aaa}) and the fact that
$\int |y| |\phi'(y)|dy = 1$. $\qed$

In Sections \ref{proof-zero} and  \ref{proof-size}, $h_t$ is given
by (\ref{def-ht}) and $f$ is the bounded solution of the Stein
equation (\ref{basic}) with $\mu=0,\sigma^2=1$, and test function
$h_t$. With $||\cdot||$ the sup norm, Lemma 3 of [\ref{Stein86}]
gives
\bea
\label{Stein-86-bounds}
||f|| \le \sqrt{2 \pi} \le 2.6
\quad
\mbox{and}
\quad ||f'|| \le 4.
\ena

\subsection{Proof of Theorem 1.1 (zero biasing)}
\label{proof-zero}
\begin{lemma}
\label{delta-independent} Let $Y$ be a mean-zero random variable
with variance $\sigma^2$, and let $Y^*$ be defined on the same
space as $Y$, with the $Y$-zero biased distribution, satisfying
$|Y^*-Y|/\sigma \le A$ for some $A$. Then
\beas
\nonumber \delta_t  \le  (6.6 +a )A+2A^2  + \frac{1}{t}\left(
2\delta A + a A^2 \right) \quad \mbox{for all $t \in (0,1)$}.
\enas
\end{lemma}

\noindent {\bf Proof:} Let $W=Y/\sigma$, whence $W^*=Y^*/\sigma$
and $|W^*-W| \le A$. By differentiation in (\ref{basic}) and
(\ref{def-ht}) respectively, we have
\bea
\label{htprime}
f''(x) = f(x) + xf'(x) + h_t'(x),
\quad
\mbox{with}
\quad
h_t'(x) = -\frac{1}{t}\int h(x+ty)\phi'(y)dy.
\ena
By (\ref{W-to-W*-via-f}) and (\ref{htprime}), with $N_th=\E
h_t(Z)$ for $Z$ a standard normal variable, we also have
\bea
\nonumber |\E  h_t(W)-Nh_t| &=& |\E  [f'(W^*)-f'(W)]| = |\E
\int_W^{W^*}f''(x)dx|
\\
\label{three-terms}
               &=& |\E  \int_W^{W^*} \left( f(x) + xf'(x) + h_t'(x) \right)
               dx|.
\ena
Let $V=W^*-W$. Applying the triangle inequality in
(\ref{three-terms}) and using (\ref{Stein-86-bounds}), for the
first term we find that
\bea
\label{first-bound} |\E  \int_W^{W^*} f(x) dx |\le 2.6 \E  |V| \le
2.6 A
\ena
and for the second term, again using (\ref{Stein-86-bounds}), and,
now,  $\E|W| \le (\E W^2)^{1/2}=1$, we find that
\bea
\nonumber && \left|\E  \int_W^{W^*} xf'(x) dx\right| \le
4\E\left|\int_W^{W+V} |x| dx\right| =2\E \left| (W+V)|W+V|-W|W|
\right| \\
&\le& \E \left( 4|WV| + 2V^2 \right) \le 4A\E |W| + 2A^2 \le 4A +
2A^2. \label{second-bound}
\ena

For the final term in (\ref{three-terms}), with $U \sim {\cal
U}[0,1]$ independent of $W$ and $V$, we write
\beas
|\E  \int_W^{W^*}h_t'(x)dx| = |\E  V \int_0^1 h_t'(W+ u V) du| =
|\E  Vh_t'(W + U V)|.
\enas
Then, using (\ref{htprime}), $\int \phi'(y)dy=0$, and Lemma
\ref{bhat-smoothing}, we have
\bea
\nonumber && |\E  Vh_t'(W + U V)| = \frac{1}{t}|\E  V \int h(W+U V +ty) \phi'(y) dy| \\
\nonumber &=& \frac{1}{t}|\E  V \int [h(W+U V +ty)-h(W+U V)] \phi'(y) dy| \\
\nonumber &\le& \frac{1}{t} \E  \left( |V| \int [h^+_{|V|
+t|y|}(W)-h^-_{|V| +t|y|}(W)]\, |\phi'(y)| dy \right)
\le \frac{1}{t}A \E  \left( \int \tilde{h}_{A +t|y|}(W) |\phi'(y)| dy \right)\\
&\le& \frac{1}{t}A \left(2\delta + a(A+t)\right) = \frac{1}{t}(2
\delta A + a A^2) + aA. \label{third-bound}
\ena
By combining the bounds (\ref{first-bound}), (\ref{second-bound}),
and (\ref{third-bound}) we complete the proof. $\qed$

\noindent {\bf Proof of Theorem \ref{delta-bound-theorem}:}
Letting $t=\alpha A$ in Lemma \ref{delta-independent}, we have
\bea
\delta_t &\le&  (6.6 +a )A+2A^2 +
                \frac{1}{\alpha A}\left( 2\delta A + a A^2
                \right)
\label{delta-t-bound}
           = (6.6 +a + \frac{a}{\alpha})A+2A^2 + \frac{2\delta}{\alpha}.
\ena
Substituting (\ref{delta-t-bound}) into the bound for $\delta $
given by Lemma \ref{bhat-smoothing}, we have
\beas
\delta &\le& 2.8((6.6 +a + \frac{a}{\alpha})A+2A^2 +
\frac{2\delta}{\alpha}) + 4.7a
\alpha A\\
&\le& 18.5 A +2.8 aA + 2.8 \frac{a A}{\alpha} +5.6 A^2 + 5.6
\frac{\delta}{\alpha} + 4.7a \alpha A,
\enas
meaning that
\bea
\label{minimize-this}
\delta \le A \left( \frac{18.5 +5.6 A +2.8 a
+ 2.8 a /\alpha + 4.7a \alpha}{1-5.6/\alpha} \right).
\ena
Setting $\alpha = 2 \times 5.6$, for which $t < 1$ since $A \le
1/12$, we obtain (\ref{delta-bound}) and, hence, the theorem.
$\qed$

\subsection{Proof of Theorem 1.2 (size biasing)}
\label{proof-size}
\begin{lemma}
Let $Y \ge 0$ be a random variable with mean $\etamu$ and variance
$\sigma^2$, and let $Y^s$ be defined on the same space as $Y$,
with the $Y$-size-biased distribution, satisfying $|Y^s-Y|/\sigma
\le A$ for some $A$. Then for all $t \in (0,1)$,
\bea
\label{lemma-size-delta-independent-bound} \delta_t  \le
\frac{\etamu}{\sigma} \left( \frac{4 \Delta}{\sigma}
+(3.3+\frac{1}{2}a)A^2+\frac{2}{3}A^3 + \frac{1}{2t}(2 \delta A^2
+ a A^3) \right),
\ena
with $\Delta$ as in (\ref{def-Delta}).
\end{lemma}

\noindent {\bf Proof:} With $W=(Y-\mu)/\sigma$, let
$W^s=(Y^s-\mu)/\sigma$ (which is a slight abuse of notation).
Then, $|W^s-W| \le A$. Note that
\bea
\label{size-bias-stein-equation-connection}
\E Wf(W) =
\frac{\etamu}{\sigma}(f(W^s)-f(W)),
\ena
and, so, with $V=W^s-W$, we have
\bea
\nonumber
&& \E  h_t(W)-Nh_t = \E  \left( f'(W) - W f(W) \right)
= \E  \left( f'(W) - \frac{\etamu}{\sigma} (f(W^s)-f(W) \right)\\
\nonumber &=& \E  \left( f'(W) - \frac{\etamu}{\sigma}
\int_W^{W^s}f'(x)dx \right) = \E  \left( f'(W) -
\frac{\etamu}{\sigma} V \int_0^1f'(W+uV)du \right)\\
\label{four-terms-in-two-groups} &=& \E  \left(  f'(W) -
\frac{\etamu}{\sigma} V f'(W) \right) + \E
\left(\frac{\etamu}{\sigma} V f'(W) -\frac{\etamu}{\sigma} V
\int_0^1f'(W+uV)du \right).
\ena
Since $ \E  (\etamu V/\sigma) = \etamu\E  (W^s-W)/\sigma = \E
(\etamu Y^s-\etamu Y)/\sigma^2 =1, $ for the first expectation in
(\ref{four-terms-in-two-groups}) we have
\bea
\E  \left\{ f'(W) \E  \left( 1 - \frac{\etamu}{\sigma} V|W\right)
\right\} \le 4 \frac{\etamu}{\sigma} \sqrt{\mbox{Var}(\E
(W^s-W|W))}=4 \frac{\etamu}{\sigma^2}\Delta,
\label{first-two-bound}
\ena
using (\ref{Stein-86-bounds}) and (\ref{def-Delta}). Now, using
(\ref{htprime}), we write the second expectation in
(\ref{four-terms-in-two-groups}) as
\bea
\nonumber && \frac{\etamu}{\sigma}  V \left\{f'(W) - \int_0^1
f'(W+uV)du \right\}  = \frac{\etamu}{\sigma}  V \int_0^1 (f'(W)-f'(W+uV))du\\
\label{three-size-bias} &=& -\frac{\etamu}{\sigma}  V \int_0^1
\int_W^{W+uV}f''(v)dv du= -\frac{\etamu}{\sigma}  V \int_0^1
\int_W^{W+uV}(f(v)+vf'(v)+h_t'(v))dv du.
\ena
We apply the triangle inequality and bound the three resulting
terms separately. For the expectation arising from the first term
on the right-hand side of (\ref{three-size-bias}), by
(\ref{Stein-86-bounds}) we have
\bea
\label{size-1of3} |\E  \{ \frac{\etamu}{\sigma}  V \int_0^1
\int_W^{W+uV} f(v) dv du\} | \le 2.6 \frac{\etamu}{\sigma}  \E
\{|V| \int_0^1 u|V| du \} \le 1.3 \frac{\etamu}{\sigma}  A^2
\ena
and, for the second term, arguing as in (\ref{second-bound}) we
have
\bea
\nonumber &&|\E  \frac{\etamu}{\sigma}  V \int_0^1 \int_W^{W+uV} v
f'(v) dv du| \le
2\frac{\etamu}{\sigma}  \E  |V| \int_0^1 \left|\int_W^{W+uV} 2 |v| dv \right| du\\
\nonumber &\le& 2\frac{\etamu}{\sigma}  \E  |V| \int_0^1
(2u|WV|+u^2V^2) du \le 2\frac{\etamu}{\sigma}  A \int_0^1 (2Au \E
|W|+u^2A^2) du \\ &\le& 2\frac{\etamu}{\sigma}  A (A + A^2/3).
\label{size-2of3}
\ena

For the last term in (\ref{three-size-bias}), the computation is
more involved than, yet similar to, that for zero biasing.
Beginning with the inner integral, we have
\beas
\int_W^{W+uV}h_t'(v)dv = u V  \int_0^1 h_t'(W+ xuV) dx
\enas
and using (\ref{htprime}),
$$\int \phi'(y)dy=0,$$
and Lemma \ref{bhat-smoothing}, for the last term in
(\ref{three-size-bias}) we have
\bea
\nonumber && |\frac{\etamu}{\sigma}  \E  \int_0^1 \int_0^1 u V^2 h_t'(W + xu V)dx du| \\
\nonumber &=& \frac{\etamu}{\sigma t}|\E  V^2\int_0^1 \int_0^1 \int u h(W+ xu V +ty) \phi'(y) dy dx du| \\
\nonumber &=& \frac{\etamu}{\sigma t}|\E  V^2 \int_0^1 \int_0^1
\int u [h(W+ xu V +ty)-h(W+ xu V)] \phi'(y) dy dx du |\\
\nonumber &\le& \frac{\etamu}{\sigma t}
\E  \left( V^2  \int \int_0^1 u [h^+_{|V| +t|y|}(W)-h^-_{|V| +t|y|}(W)]\, |\phi'(y)| du dy \right) \\
\nonumber &=& \frac{\etamu}{2\sigma t}
\E  \left( V^2  \int [h^+_{|V| +t|y|}(W)-h^-_{|V| +t|y|}(W)]\, |\phi'(y)| dy \right) \\
&\le& \nonumber
\frac{\etamu}{2\sigma t}A^2  \E  \left( \int \tilde{h}_{A +t|y|}(W) |\phi'(y)| dy \right)\\
\nonumber &\le& \frac{\etamu}{2\sigma t}A^2 \left(2\delta
+ a(A+t)\right) \\
\label{size-3of3} &=& \frac{\etamu}{2\sigma t}(2\delta A^2 + a
A^3) + \frac{\etamu}{2\sigma}aA^2.
\ena
By combining (\ref{first-two-bound}), (\ref{size-1of3}),
(\ref{size-2of3}), and (\ref{size-3of3}) we complete the proof.
$\qed$

\noindent {\bf Proof of Theorem \ref{size-delta-bound-theorem}}
Applying Lemma \ref{bhat-smoothing} using the bound
(\ref{lemma-size-delta-independent-bound}) on $\delta_t$, we have
$$
\delta \le 2.8 \frac{\etamu}{\sigma} \left( \frac{4\Delta}{\sigma}
+(3.3+\frac{1}{2}a)A^2+\frac{2}{3}A^3 + \frac{1}{2t}(2\delta A^2 +
a A^3) \right) + 4.7 at,
$$
or,
\bea
\label{minimize-this-size} \delta \le \frac{2.8
(\etamu/\sigma)\left( 4\Delta/\sigma
+(3.3+\frac{1}{2}a)A^2+\frac{2}{3}A^3 + a A^3/2t \right) + 4.7
at}{1-2.8 \etamu A^2 /(\sigma t)}.
\ena
Setting $t=2 \times 2.8 \etamu A^2/\sigma$, such that $t <1$ since
$A \le (\sigma/(6 \etamu))^{1/2}$, (\ref{size-delta-bound}) now
follows from
\beas
\delta &\le& 5.6 \frac{\etamu}{\sigma}\left(
\frac{4\Delta}{\sigma} +(3.3+\frac{1}{2}a)A^2+\frac{2}{3}A^3 +
\frac{\sigma }{2(5.6 \etamu)}a A
\right) + 2(4.7)a (5.6 \frac{\etamu A^2}{\sigma})\\
&\le& \frac{aA}{2} + \frac{\etamu}{\sigma} \left((19+56a)A^2 +
4A^3 \right)+ 23 \frac{\etamu \Delta}{\sigma^2}. \,\,\,\,\,\qed
\enas

There are compromises in the choice of smoothing parameter; if we
take $\alpha = 4 \times 5.6$ in (\ref{minimize-this}) for $B \le
\sigma/48$, and $t=4 \times 2.8 \etamu A^2/\sigma$ in
(\ref{minimize-this-size}) for $B \le \sigma^{3/2}/(12
\etamu)^{1/2}$, bounds (\ref{delta-bound}) and
(\ref{size-delta-bound}) become
\bea
\label{delta-bound-alt} \delta &\le& A(145 a + 7.5A + 25)
\ena
and
\bea
\label{size-delta-bound-alt} \delta &\le& \frac{aA}{6} +
\frac{\etamu}{\sigma} \left((13+73a)A^2 + 2.5A^3 \right)+15
\frac{\mu \Delta}{\sigma^2},
\ena
respectively.

\section{Remarks}
The zero- and size-bias coupling both conform well to Stein's
characterizing equation, and their use produces bounds on the
distance of a random variable $Y$ to the normal in many instances.
The couplings are adaptable to the situation; in particular, the
size-biased coupling, previously used in [\ref{Goldstein-Rinott}]
for global dependence, is applied here to handle cases of local
dependence.

The applications in Section \ref{zero-applications} illustrate how
bounds on the distance $\delta$ from $Y$ to the normal can be
generated using only a zero-bias coupling and a bound on
$|Y^*-Y|$; in particular, the bounds do not depend on the
often-difficult calculation of variances of conditional
expectations of the form $\mbox{Var}\{\E ({\tilde Y}-Y|Y)\}$,
which appear in the exchangeable-pair and size-biased versions of
Stein's method when coupling $Y$ to some ${\tilde Y}$. It is hoped
that this feature of the zero-bias method will motivate a better
understanding of the construction of couplings of $Y^*$ to $Y$ in
greater generality than those that depend on the existence of the
exchangeable pair of Proposition \ref{peara}. In particular, the
applications in Section \ref{size-applications} show an evidently
wider scope of applicability of the size bias coupling over the
zero bias one, as it is presently understood.

\begin{center}
{\bf Acknowledgments}
\end{center}
The author would like to thank Martin Raic and Qi-Man Shao for
their insightful comments on an earlier version of this
manuscript.

\section*{Bibliography}
\begin{enumerate}

\item \label{BRS} Baldi, P., Rinott, Y., and Stein, C. (1989) A
normal approximation for the number of local maxima of a random
function on a graph. Probability, statistics, and mathematics,
59-81, Academic Press, Boston, MA.

\item \label{BR} Baldi, P. and Rinott, Y. (1989) Asymptotic normality of some
graph-related statistics.  {\em J. Appl. Probab.}  {\bf 26},
171-175.

\item \label{Barbour} Barbour, A. D., Janson, S.,  Karo\'nski,
M., and Ruci\'nski, A. (1990) Small cliques in random graphs. {\em
Random Structures and Algorithms}, {\bf 1}, 403-434.

\item \label{Dia} Bayer, D. and Diaconis, P. (1992) Trailing the Dovetail Shuffle to its
Lair. {\em The Annals of Applied Probability}, {\bf 2}, 294-313.

\item \label{bhat} Bhattacharya, R.N. and Ranga Rao R. (1986). {\em Normal
approximation and asymptotic expansion}, Krieger, Melbourne, Fla.

\item \label{Biggs} Biggs, N. (1993). {\sl Algebraic Graph Theory},
Cambridge University Press.

\item \label{Bolt} Bolthausen, E. (1984) An estimate of the reminder
in a combinatorial central limit theorem. {\em Z. Wahrsch. Verw.
Gebiete.}, {\bf 66}, 379-386.

\item \label{BG93} Bolthausen, E. and G\"otze, F. (1993).
The rate of convergence for multivariate sampling statistics. {\em
Ann. Statist.} {\bf 21}, 1692-1710.

\item \label{Brouwer} Brouwer, A.E., Cohen, A. M.,  and Neumaier, A. (1989).
{\em Distance-Regular  Graphs}, Springer-Verlag, Berlin.

\item \label{Chen-Shao} Chen, Louis H.Y., and Shao, Q.M. (2004)
Normal approximation under local dependence. {\em Ann Prob.} {\bf
32}, 1985-2028.

\item \label{DW-1} Darling, R. W. R., and Waterman, M.S. (1986) Extreme value
distribution for the largest cube in a random lattice.  {\em SIAM
J. Appl. Math.}  {\bf 46}, 118-132.

\item \label{DW-2} Darling, R. W. R., and Waterman, M.S. (1985) Matching rectangles in
$d$ dimensions: algorithms and laws of large numbers. {\em Adv. in
Math.} {\bf 55}, 1-12.

\ignore{Dembo, Amir; Rinott, Yosef Some examples of normal
approximations by Stein's method.  Random discrete structures
(Minneapolis, MN, 1993),  25--44, IMA Vol. Math. Appl., 76,
Springer, New York, 1996.}



\item \label{naus-book} Glaz, J.,  Naus, J., and Wallenstein, S.
(2001) {\em Scan statistics.} Springer Series in Statistics.
Springer-Verlag, New York.

\item \label{res} Goldstein, L. (2004) Normal Approximation for
Hierarchical Sequences, {\em Annals of Applied Probability}, {\bf
14}, pp. 1950-1969. arXiv:math.PR/0503549

\item \label{Goldstein-Reinert} Goldstein, L. and Reinert, G.
(1997) Stein's Method and the Zero Bias Transformation with
Application to Simple Random Sampling, {\em Annals of Applied
Probability}, {\bf 7}, 935-952. arXiv:math.PR/0510619


\item \label{Goldstein-Reinert-03}  Goldstein, L. and Reinert, G.
(2005) Distributional transformations, orthogonal polynomials, and
Stein characterizations, {\em Journal of Theoretical Probability},
{\bf 18}, 185-208. arXiv:math.PR/0510240

\item \label{Goldstein-Rinott} Goldstein, L. and Rinott, Y.
(1996). Multivariate normal approximations by Stein's method and
size bias couplings, {\em J. Appl. Prob.} {\bf 33}, 1-17.
arXiv:math.PR/0510586

\item \label{Goldstein-Rinott-Perm} Goldstein, L. and Rinott, Y.
(2004). A permutation test for matching and its asymptotic
distribution, {\em Metron}, {\bf 61} (2003), pp. 375-388.
arXiv:math.PR/0510240



\item \label{Gotze} G\"otze, F. (1991).
On the rate of convergence in the multivariate CLT. {\em Annals of
Probability}, {\bf 19}, 724-739.

\item \label{ChenHo} Ho, S. T. and Chen, Louis H. Y. (1978) An $L\sb{p}$
bound for the remainder in a combinatorial central limit theorem.
{\em Ann. Probab.} {\bf 6}, 231-249.

\item \label{huang} Huang, H. (2002) Error bounds on multivariate normal
approximations for word count statistics. {\em Adv. in Appl.
Probab.} {\bf 34}, 559-586.

\item \label{Kolchin} Kolchin, V.F., and Chistyakov, V.P. (1973) On a
combinatorial limit theorem. {\em Theory Probability Appl.}, {\bf
18}, 728-739.

\item \label{naus-scan} Naus, J. I. (1982)
Approximations for distributions of scan statistics. {\em J. Amer.
Statist. Assoc.}  {\bf 77}, 177-183.

\item \label{Rin-Rot-1} Rinott, Y. and Rotar, V. (1996). A multivariate CLT for local dependence with
$n\sp {-1/2}\log n$ rate and applications to multivariate graph
related statistics. {\em J. Multivariate Anal.},  {\bf 56},
333-350.

\item \label{Rin-Rot-2} Rinott, Y. and Rotar, V. (1997). On coupling constructions and rates in the CLT for dependent
summands with applications to the antivoter model and weighted
$U$-statistics. {\em Ann. Appl. Probab.},  {\bf 7}, 1080-1105.

\item \label{Stein72} Stein, C. (1972). A bound for the error in the normal
approximation to the distribution of a sum of dependent random
variables. Proc. Sixth Berkeley Symp. Math. Statist. Probab. {\bf
2}, 583-602, Univ. California Press, Berkeley.

\item \label{Stein81} Stein, C. (1981). Estimation of the mean of
a multivariate normal distribution. {\em Ann. Statist.}, {\bf 9},
1135-1151.

\item \label{Stein86} Stein, C.  (1986). {\em Approximate
Computation of Expectations.} IMS, Hayward, CA.

\item \label{Bahr} von Bahr, B. (1976) Remainder term estimate in
a combinatorial limit theorem. {\em Z. Wahrsch. Verw. Gebiete.},
{\bf 35}, 131-139.

\end{enumerate}
\end{document}